\documentclass[10pt]{amsart}
\usepackage[all]{xy}
\usepackage{amsmath}
\usepackage{amsfonts}
\usepackage{amssymb}
\usepackage{amscd}
\usepackage{amsthm}
\usepackage{latexsym}
\usepackage{amsbsy}
\usepackage{graphicx}
\usepackage{epstopdf}
\AppendGraphicsExtensions{.gif}
\newtheorem{lemma}{Lemma}[section]
\newtheorem{proposition}{Proposition}[section]
\newtheorem{theorem}{Theorem}[section]
\newtheorem{corollary}{Corollary}[section]
\newtheorem{definition}{Definition}[section]

\newtheorem{remark}{Remark}[section]

\begin{document}

\title{Scalar Curvature for Noncommutative Four-Tori}

\author{FARZAD FATHIZADEH AND MASOUD KHALKHALI}

\date{}

\maketitle

\begin{abstract}

In this paper we study the curved geometry of noncommutative $4$-tori 
$\mathbb{T}_\theta^4$. We use a Weyl conformal factor to perturb the standard volume 
form and obtain the Laplacian that encodes the local geometric information. 
We use Connes' pseudodifferential calculus to explicitly compute the 
terms in the small time heat kernel expansion of the perturbed Laplacian 
which correspond to the volume and scalar curvature of  $\mathbb{T}_\theta^4$. 
We establish the analogue of Weyl's law, define a noncommutative residue, 
prove the analogue of Connes' trace theorem,  and find explicit 
formulas for the local functions that describe the scalar 
curvature of $\mathbb{T}_\theta^4$. We also study the analogue of the Einstein-Hilbert 
action for these spaces and show that metrics with constant scalar curvature 
are critical for this action. 
$ $\\
$ $\\ 

\end{abstract}

\tableofcontents

\section{Introduction}

Spectral geometry has played an important role in the development of metric 
aspects of noncommutative geometry \cite{con1.5, con3, conmos1, conmar}.  
After the seminal paper \cite{contre}, in which the analogue of the Gauss-Bonnet 
theorem is proved for noncommutative two tori $\mathbb{T}_\theta^2$, there 
has been much progress in understanding the local differential geometry of 
these noncommutative spaces \cite{fatkha1, conmos2, fatkha2, fatkha3}. In these works, the  
flat geometry  of $\mathbb{T}_\theta^2$ which was studied in \cite{con1}  is conformally 
perturbed by means of a Weyl factor given by a positive invertible 
element in $C^\infty(\mathbb{T}_\theta^2)$ (see also \cite{cohcon} for a preliminary 
version). Connes' 
pseudodifferential calculus developed in \cite{con1} for $C^*$-dynamical systems 
is employed crucially to apply heat kernel techniques to geometric operators 
on $\mathbb{T}_\theta^2$ to derive small time heat kernel expansions that 
encode local geometric information such as scalar curvature.  A purely 
noncommutative feature is the appearance of the modular automorphism of 
the state implementing the conformal perturbation of the metric 
in the computations and in the final formula for the curvature 
\cite{conmos2, fatkha2}.

In this paper we study the curved geometry of noncommutative $4$-tori 
$\mathbb{T}_\theta^4$. We view these spaces as noncommutative abelian 
varieties equipped with a complex structure given by the simplest possible period matrix.  
We use a Weyl conformal factor to perturb the standard volume 
form and obtain the Laplacian that encodes the local geometric information. 
We use Connes' pseudodifferential calculus to explicitly compute the 
terms in the small time heat kernel expansion of the perturbed Laplacian 
which correspond to the volume and scalar curvature of  $\mathbb{T}_\theta^4$. 
We establish the analogue of Weyl's law, define a noncommutative residue, 
prove the analogue of Connes' trace theorem \cite{con2},  and find explicit 
formulas for the local functions that describe the scalar 
curvature of $\mathbb{T}_\theta^4$. We also study the analogue of the 
Einstein-Hilbert action for these spaces and show that metrics with constant 
scalar curvature are critical for this action.

This paper is organized as follows. In Section \ref{NCT}, we recall basic facts about higher 
dimensional noncommutative tori and their flat geometry. In Section \ref{LHK}, 
we consider the noncommutative 4-torus $\mathbb{T}_\theta^4$ with the simplest 
structure of a noncommutative abelian variety. We perturb the standard 
volume form on this space conformally (cf. \cite{contre}), and analyse the 
corresponding perturbed Laplacian. Then, we recall Connes' pseudodifferential 
calculus \cite{con1} for $\mathbb{T}_\theta^4$ and review the derivation of  the 
small time heat kernel expansion for the perturbed Laplacian, using this calculus. 
In Section \ref{WLCTT}, we prove the analogue of Weyl's law for $\mathbb{T}_\theta^4$ 
by studying the asymptotic distribution of the eigenvalues of the perturbed Laplacian on 
this space. We then define a noncommutative residue on the algebra of classical 
pseudodifferential operators  on $\mathbb{T}_\theta^4$, and show 
that it gives the unique continuous trace on this algebra. We also prove the analogue of 
Connes' trace theorem for $\mathbb{T}_\theta^4$ by showing that this noncommutative 
residue and the Dixmier trace coincide on pseudodifferential operators of order $-4$. 
In Section \ref{SCEHA}, we perform the computation of the scalar curvature for 
$\mathbb{T}_\theta^4$, and find explicit formulas for the local functions that describe 
the curvature in terms of the modular automorphism of the conformally perturbed volume 
form and derivatives of the logarithm of the Weyl factor. Then, by integrating this curvature, 
we define and find an explicit formula  for the analogue of the Einstein-Hilbert  action for 
$\mathbb{T}_\theta^4$. Finally, we show that the extremum of this action occurs at  
metrics with constant scalar curvature (see \cite{braors} for the corresponding 
commutative statement).

We are indebted to Alain Connes for several enlightening and productive 
discussions at different stages of this work. Also, F. F. would like to thank 
IHES for the excellent environment and kind support during his visit in 
Winter 2012, where part of this work was carried out.

\section{Noncommutative Tori} \label{NCT}

In this section we recall basic facts about higher 
dimensional noncommutative tori and their flat geometry.

\subsection{Noncommutative real tori.} \label{ncrtori}

Let $V$ be a finite dimensional real vector space equipped 
with a positive definite inner product $\langle \, , \, \rangle$ and 
let $\theta: V \otimes V \to \mathbb{R}$ be a  skew-symmetric 
bilinear  form on $V$. Let $\Lambda \subset V$ be a cocompact  
lattice in $V$. Thus $\Lambda$  is a discrete abelian subgroup 
of $V$  such that  the quotient space $V/ \Lambda$ is compact.  
Equivalently, we can describe 
$\Lambda$ as  $\Lambda = \mathbb{Z}e_1+  \cdots  +\mathbb{Z} e_n$, 
the free abelian group generated by a linear  
basis  $e_1, \dots, e_n$ for $V$.

By definition, the noncommutative torus $C(\mathbb{T}_\theta^n)$, 
attached to the above data,  is the universal unital $C^*$-algebra 
generated by unitaries $U_{\alpha}, \alpha \in \Lambda$, satisfying the 
relations
$$ 
U_{\alpha}U_{\beta}= e^{ \pi  i \theta (\alpha, \beta)} U_{\alpha + \beta}, 
\qquad  \alpha, \beta \in \Lambda.
$$
Let  $e_1, \dots, e_n$  be a basis for $\Lambda$ over $\mathbb{Z}$, and 
let $U_i=U_{e_i}$. Then we have  
$$ 
U_k U_l = e^{2 \pi  i \theta_{ kl}}U_lU_k, \qquad  k, l=1,\dots, n,
$$
where $\theta_{ kl}=\theta (e_k, e_l)$.

Let $\Lambda' \subset V$ denote the dual lattice. Thus 
$v \in  \Lambda'$ iff $ \langle v , w \rangle \in 
2 \pi \mathbb{Z}$ for all $w \in \Lambda$. There is a continuous action 
of the dual  torus $V  / \Lambda' $  on $C(\mathbb{T}_\theta^n)$  by $C^*$-algebra 
automorphisms $\{\lambda_s \}_{s \in V}$,  defined by  
$$ 
\lambda_s  (U_{\alpha})= e^{ i \langle s, \alpha \rangle }U_{\alpha}. 
$$ 
The space of smooth elements  of this action, namely those elements 
$a \in C(\mathbb{T}_\theta^n)$ for which the map $ s \to \lambda_s (a)$ is smooth  
will be denoted by $C^\infty(\mathbb{T}_\theta^n)$.  It is a dense $*$-subalgebra of 
$C(\mathbb{T}_\theta^n)$ which plays the role of smooth functions on the noncommutative 
torus $\mathbb{T}_\theta^n$. It can be alternatively described as the algebra of elements 
in $C(\mathbb{T}_\theta^n)$ whose (noncommutative) Fourier expansion has rapidly 
decreasing Schwartz class coefficients:
\[
C^\infty(\mathbb{T}_\theta^n)=\Big \{\sum_{\alpha \in \Lambda }a_{\alpha}U_{\alpha}; 
\qquad \sup_{\alpha \in \Lambda} (|\alpha|^k  |a_{\alpha}|) < \infty, \,\, 
\forall k  \in  \mathbb{N}  \Big \}.
\]
There is a  normalized faithful  positive trace,  i.e. a tracial state,    
$\varphi_0$  on $C(\mathbb{T}_\theta^n)$ whose restriction on smooth elements is given by
\[
\varphi_0 \Big (\sum_{\alpha \in \mathbb{Z}^n}a_{\alpha}U_{\alpha} \Big )=a_{0}.
\]

The infinitesimal generator of the action  $\lambda_s$   defines a Lie algebra 
map
\begin{equation} \label{infinitesimal}
\delta : V \to \text{Der}(C^\infty(\mathbb{T}_\theta^n), C^\infty(\mathbb{T}_\theta^n)),
\end{equation}
where we have used the natural identification of the (abelian) Lie algebra 
of the torus  $V  / \Lambda'$ with $V$.

Let us fix an orthonormal basis $e_1, \dots, e_n$ for $V$. Then the 
restriction of the above map $\delta$ defines commuting derivations
$\delta_i := \delta (e_i)  : C^\infty(\mathbb{T}_\theta^n) \to C^\infty(\mathbb{T}_\theta^n),  i=1, 
\dots, n, $ which satisfy 
\[
\delta_i(U_j)= \delta_{ij}U_i, \qquad     i,j=1, \dots, n.
\]
The derivations  $\delta_j$ are analogues of the differential operators
$\frac{1}{i}\partial/\partial x_j$ acting on  smooth functions on
the ordinary  torus. We have $\delta_j(a^*)= -\delta_j(a)^* $ for $j=1, \dots, n,$ 
and any $a\in C^\infty(\mathbb{T}_\theta^n)$. Moreover, since $\varphi_0 \circ \delta_j =0$, 
for all  $j$, we have the analogue of the integration by parts formula:
\[ 
\varphi_0(a\delta_j(b)) = -\varphi_0(\delta_j(a)b), \qquad  a,b \in 
C^\infty(\mathbb{T}_\theta^n). 
\]

Using these derivations, we can define the flat Laplacian 
$$ 
\triangle= \sum_{i=1}^n \delta_i^2: 
C^\infty(\mathbb{T}_\theta^n) \to C^\infty(\mathbb{T}_\theta^n).
$$
We note that the Laplacian $\triangle$ is independent of the choice of 
the orthonormal basis $e_1, \dots, e_n.$

\subsection{Noncommutative complex tori.}

Let  $W$ be an $n$-dimensional complex vector space and 
$\Lambda \subset W$ be a lattice in $W$. Thus $\Lambda$ is a 
free  abelian group of rank $2n$ which is discrete  in its subspace 
topology.  Given a basis $e_1, \dots, e_n$ of $W$ as a complex 
vector space and a basis $\lambda_1, \dots, \lambda_{2n}$ of 
$\Lambda$ as a free abelian group, we can express 
$\lambda_1, \dots, \lambda_{2n}$ in terms of $e_1, \dots, e_n$. 
We obtain an  $n$ by $2n$ matrix $ \mathcal{M}= (A, B)$ with 
$A, B \in M_n (\mathbb{C})$ with 
$$
\lambda_j= \sum_{i=1}^n  \mathcal{M}_{ij}e_i, \qquad  j=1, \dots, 2n.
$$

Let $W_{\mathbb{R}}$ denote the realification of $W$. Note 
that $\lambda_1, \dots, \lambda_{2n}$ is a basis for $W_{\mathbb{R}}$ 
over $\mathbb{R}$.  Let $dz_1, \dots, dz_n$ denote the  basis of 
$W^*= \text{Hom}_{\mathbb{C}}\, (W, \mathbb{C})$,  
dual to the basis  $e_1, \dots, e_n$, and let 
$dx_1, \dots, dx_{2n}$ denote the basis of 
$W_{\mathbb{R}}^* =  \text{Hom}_{\mathbb{R}}\, (W, \mathbb{R}),$ 
dual  to the basis  $\lambda_1, \dots, \lambda_{2n}$. Then we have for 
$i=1, \dots, n$,  
$$ 
dz_{i}= \sum_{j=1}^{2n} \mathcal{M}_{ij}dx_j, \qquad   
d \bar{z}_{i}= \sum_{j=1}^{2n}  \bar{\mathcal{M}}_{ij}dx_j. 
$$

Now let $\theta :W_{\mathbb{R}} \otimes W_{\mathbb{R}} \to \mathbb{R}$ 
be an alternating bilinear form and  $\mathbb{T}_{\theta}^n$ 
denote the associated noncommutative torus.  We further assume  that 
$W_{\mathbb{R}}$ is equipped with an Euclidean inner product. Using \eqref{infinitesimal}, 
we get  the derivations
$$
\delta_i = \delta_{\lambda_i}: C^\infty(\mathbb{T}_\theta^n)\to 
C^\infty(\mathbb{T}_\theta^n), \qquad  i=1, \dots, 2n.   
$$

The above relations define the Dolbeault operators which will be denoted 
by $\partial_i, \bar \partial_i, i=1, \dots, n.$

We have a  decomposition   
$W_{\mathbb{R}} \otimes \mathbb{C}= W_{(1, 0)} \oplus W_{(0,1)}$  
with  $W_{(1, 0)} =W$ and $W_{(0, 1)} = \bar{W}$ with  the bases   
$e_1, \dots, e_n$  for $W_{(1, 0)}$ and $ \bar{e}_1, 
\dots,  \bar{e}_n$ for $W_{(0,1)}.$ Let $dz_i$ and $d\bar{z}_j$ 
denote the corresponding bases for  dual spaces  
$W_{(1, 0)}^*= \text{Hom}_{\mathbb{C}} (W_{(1, 0)}, \mathbb{C})$ 
and  $W_{(0, 1)}^*.$ Using this decomposition we can define a  
Dolbeault type complex for on  $C^\infty(\mathbb{T}_\theta^n)$  as follows. Let
$$
\Omega^{p, q}: = C^\infty(\mathbb{T}_\theta^n)\otimes \wedge^pW_{(1, 0)}^*\otimes 
\wedge^qW_{(0, 1)}^*, 
$$
and define the operators
$$ 
\partial_i: \Omega^{p, q} \to \Omega^{p+1, q}, \qquad     
\bar{ \partial}_i: \Omega^{p, q} \to \Omega^{p, q+1}.
$$
by 
\begin{eqnarray*}
\partial_i (a dz_I \wedge d\bar{z}_J) &=&\sum_i \partial_i (a) dz_i 
\wedge dz_I \wedge d\bar{z}_J,  \\
\bar{\partial}_i (a dz_I \wedge d\bar{z}_J) &=&\sum_i  
\bar{\partial}_i (a) d \bar{z}_i \wedge dz_I \wedge d\bar{z}_J.
\end{eqnarray*}
Suppressing the obvious indexing, these operators satisfy the 
relations
$$ 
\partial^2=0, \qquad  \bar{\partial}^2=0, \qquad  {\partial}\bar{\partial}+   
 \bar{\partial} {\partial}=0.
$$

Let $\frak{H}_n \subset M_n(\mathbb{C})$ denote the Siegel 
upper half space. By definition, a matrix $\Omega \in \frak{H}_n$ if 
and only if 
$$
\Omega^ t =\Omega  \qquad \text{and} \qquad \text{Im}\, \Omega >0.
$$ 
For $n=1$, $ \frak{H}_1$ is the Poincar\'e upper half plane. 
The following two conditions are known to be equivalent 
for a lattice $\Lambda \subset W$: i) The complex torus 
$W/\Lambda$  can be embedded, as a 
complex manifold, in a complex projective space 
$\mathbb{P}^N (\mathbb{C})$; ii) There exists a basis 
$(e_1, \dots, e_n)$  of $W$, and a basis 
$(\lambda_1, \dots, \lambda_{2n})$ of $\Lambda$ such that 
the matrix of  $(\lambda_1, \dots, \lambda_{2n})$ in the basis 
$(e_1, \dots, e_n)$ is of the form
$$
(\Delta_{k}, \Omega),
$$
where  $k = (k_1, \dots, k_n)$ is a sequence of integers 
$k_i \in \mathbb{Z}$, $\Delta_k = \text{diag}\, (k_1, \dots, 
k_n)$ is a diagonal matrix, and $\Omega \in \frak{H}_n$. 
A noncommutative torus $\mathbb{T}_\theta^n$ attached to a pair $(W, \Lambda)$ satisfying the equivalent 
conditions i) or ii) can be regarded as a noncommutative  abelian variety.

\section{Laplacian and its  Heat Kernel} \label{LHK}

In this section we consider the curved geometry of noncommutative 
$4$-tori and analyse the corresponding Laplacian. We also recall 
Connes' pseudodifferential calculus \cite{con1} for these spaces 
and review the derivation of small time heat kernel expansions 
using this calculus.

\subsection{Perturbed Laplacian on $\mathbb{T}_\theta^4$.} 
\label{perturbedlaplacian}

We consider the construction  of Section \ref{NCT} for the 
noncommutative 4-torus $\mathbb{T}_\theta^4$, the matrix $\Omega=i I_{2 \times 2}$ 
%$\Omega = \left(\begin{array}{cc}
%i & 0   \\
%0 & i 
%\end{array}\right)$ 
in the Siegel upper  half space $\mathfrak{H}_2$, and $k_1=k_2=1$. 
That is, we  consider the complex structure on $\mathbb{T}_\theta^4$  
which is introduced by setting 
$\partial, \bar \partial : C^\infty(\mathbb{T}_\theta^4) \to 
C^\infty(\mathbb{T}_\theta^4) \oplus C^\infty(\mathbb{T}_\theta^4)$ as 
\[
\partial  = \partial_1 \oplus \partial_2, \qquad 
\bar \partial  = {\bar \partial}_1 \oplus  {\bar \partial}_2, 
\]
where
\[ 
\partial_1 =  \delta_1 - i \delta_3, \qquad 
\partial_2=  \delta_2 - i \delta_4,
\]
\[ \bar \partial_1 =  \delta_1 + i \delta_3, 
\qquad \bar \partial_2= \delta_2 + i \delta_4.
\]
We consider the inner product  
\begin{equation} \label{innerp-vol}
(a, b) = \varphi_0(b^*a), \qquad a,b \in C(\mathbb{T}_\theta^4),
\end{equation}
and denote the Hilbert space completion of $C(\mathbb{T}_\theta^4)$ 
with respect to this inner product by $\mathcal{H}_0$. 
Since the derivations $\delta_i$ are formally selfadjoint with respect 
to the above inner product, we have  
\[
\partial^* \partial = \delta_1^2+ \delta_2^2+\delta_3^2 +\delta_4^2, 
\] 
which is the flat Laplacian $\triangle$ introduced in Section \ref{NCT}. The 
reason for this coincidence is that  the underlying metric is K\"ahler 
when we have the non-perturbed standard volume form $\varphi_0$. Therefore,  the 
ordinary Laplacian and the Dolbeault Laplacian agree with each other 
in this case.

In order to perturb the above Laplacian conformally, following 
\cite{cohcon, contre}, we consider a selfadjoint  element 
$h \in C^\infty(\mathbb{T}_\theta^4)$ and perturb the volume 
form $\varphi_0$ by replacing it with 
the linear functional 
$\varphi : C(\mathbb{T}_\theta^4) \to \mathbb{C}$ defined by 
\begin{equation} 
\varphi(a)=\varphi_0(a e^{-2h}), \qquad a \in C(\mathbb{T}_\theta^4). \nonumber
\end{equation}
This is a non-tracial state which is a twisted trace, and satisfies the 
KMS condition at $\beta = 1$ for the 1-parameter group 
$\{\sigma_t \}_{t \in \mathbb{R}}$ of inner automorphisms
\[
\sigma_t (a) = e^{2ith} a e^{-2ith}, \qquad a \in C(\mathbb{T}_\theta^4).
\]
The inner product  associated with this 
linear functional is given by 
\begin{equation} 
(a, b)_\varphi = \varphi(b^*a), \qquad a,b \in C(\mathbb{T}_\theta^4). \nonumber
\end{equation}
We denote the Hilbert space completion of $C(\mathbb{T}_\theta^4)$ 
with respect to this inner product by $\mathcal{H}_\varphi$, and 
write the analogue  of the de Rham differential  on 
$\mathbb{T}_\theta^4$ as 
\[
d = \partial \oplus \bar \partial : \mathcal{H}_\varphi \to 
\mathcal{H}_\varphi^{(1,0)} \oplus \mathcal{H}_\varphi^{(0,1)}. 
\]
Here, we view $d$ as an operator from $\mathcal{H}_\varphi$ to the 
corresponding Hilbert space completion of the analogue of 
1-forms, namely the direct sum of the linear span of  
$\{a \partial b; a, b \in C^\infty(\mathbb{T}_\theta^4) \}$ and 
$\{a \bar \partial b; a, b \in C^\infty(\mathbb{T}_\theta^4) \}$.

Now we define the perturbed Laplacian 
\[
\triangle_{\varphi}=d^*d,
\]
which is an unbounded operator acting in $\mathcal{H}_\varphi$.  We will 
see in the following lemma that $\triangle_{\varphi}$ is anti-unitarily 
equivalent to a differential operator acting in $\mathcal{H}_0$,  
and because of this equivalence, we identify these operators with 
each other in the sequel.

\begin{lemma}

The perturbed Laplacian $\triangle_{\varphi}$ 
is anti-unitarily equivalent to the operator 
\[ 
e^{h} \bar \partial_1 e^{-h} \partial_1 e^{h} + 
e^{h} \partial_1 e^{-h} \bar \partial_1 e^{h} + 
e^{h} \bar \partial_2 e^{-h} \partial_2 e^{h} + 
e^{h}  \partial_2 e^{-h} \bar \partial_2 e^{h},
\]
acting in $\mathcal{H}_0$. 

\begin{proof}

It follows easily from the argument given in the 
proof of the following lemma. 

\end{proof}

\end{lemma}

Let $h', h'' \in C^\infty( \mathbb{T}_\theta^4)$ be selfadjoint elements, 
and $H$, $H_0$, $H_1$ respectively be the Hilbert space 
completion of $C(\mathbb{T}_\theta^4)$ with respect to the  inner 
products defined by
\[
 (a, b) = \varphi_0(b^*a), 
\qquad (a, b)_0 = \varphi_0(b^*a e^{-h'}), 
\qquad (a, b)_1=  \varphi_0(b^*a e^{-h''}), 
\]
for any $a, b \in C(\mathbb{T}_\theta^4)$. 
We recall the operator 
$\partial_1=\delta_1-i\delta_3: H \to H$ and its 
adjoint $\bar \partial_1 =\delta_1+i\delta_3$.

\begin{lemma} \label{antiunit} 

Let $\partial_{0,1}$  be the same map as $\partial_1$ viewed as 
an operator from $H_0$ to $H_1$. Then its adjoint is 
given by 
\[
\partial_{0,1}^* (y) = \bar \partial_1 (y e^{-h'' }) e^{h'},
\]
and the operator $\partial_{0,1}^*\partial_{0,1} : H_0 \to H_0$ is 
anti-unitarily equivalent to 
\[
e^{h'/2}\partial_1 e^{-h''}\bar \partial_1 e^{h'/2}: H \to H,
\]
where $e^{h'}, e^{-h''}, e^{h'/2}$ act by left multiplication.

\begin{proof} 

Right multiplication by $e^{h'/2}$ extends to a unitary 
map from $H$ to $H_0$, which we 
denote by $W_0$. Similarly, right multiplication by $e^{-h''/2}$ 
extends to a unitary map from $H_1$ 
to $H$, which will be denoted  by $W_1$. So we have 
\[ 
W_1 \partial_{0,1} W_0 = R_{e^{-h''/2}} \partial_1 R_{e^{h'/2}}.
\] 
Therefore, we have 
\[  
W_0^* \partial_{0,1}^* W_1^*W_1 \partial_{0,1} W_0 = R_{e^{h'/2}} 
\bar \partial_1 
R_{e^{-h''/2}} R_{e^{-h''/2}} \partial_1 R_{e^{h'/2}}. 
\]
Thus  $\partial_{0,1}^*\partial_{0,1}$ is unitarily equivalent to 
\[
R_{e^{h'/2}} \bar \partial_1 R_{e^{-h''/2}} R_{e^{-h''/2}} \partial_1 R_{e^{h'/2}}.
\]
Conjugating the latter with the anti-unitary involution $J(a)=a^*$, one 
can see that it is anti-unitarily equivalent to 
\[
e^{h'/2}\partial_1 e^{-h''}\bar \partial_1 e^{h'/2}.
\]

\end{proof}

\end{lemma}

\subsection{Connes' pseudodifferential calculus for $\mathbb{T}_\theta^4$.}
\label{ncpseudos}

A pseudodifferential calculus was developed in \cite{con1} for 
$C^*$-dynamical systems. Here we briefly recall  this calculus for 
the canonical dynamical system defining the noncommutative 4-torus, 
and will use it in the sequel to apply heat kernel techniques 
\cite{gil, contre} to the perturbed Laplacian $\triangle_\varphi$ on 
$\mathbb{T}_\theta^4$.

A differential operator of order $m \in \mathbb{Z}_{\geq 0}$ on 
$\mathbb{T}_\theta^4$ is an operator of the form 
\[
\sum_{|\ell| \leq m} a_\ell \delta_1^{\ell_1} 
\delta_2^{\ell_2} \delta_3^{\ell_3} \delta_4^{\ell_4},  
\]
where 
$ \ell= (\ell_1, \ell_2, \ell_3, \ell_4) \in \mathbb{Z}_{\geq 0}^4,  
|\ell |=\ell_1+\ell_2+\ell_3+ \ell_4,  
a_\ell \in C^\infty(\mathbb{T}_\theta^4).$ We first recall the definition of the 
operator valued symbols, using which, the notion of differential operators on 
$\mathbb{T}_\theta^4$ extends to the notion of pseudodifferential 
operators \cite{con1}. For convenience, we will use the notation 
$\partial_j$ for the partial derivatives $\frac{\partial}{\partial \xi_j}$ 
with respect to the coordinates $\xi=(\xi_1, \dots, \xi_4) \in \mathbb{R}^4$, 
and for any $\ell \in \mathbb{Z}_{\geq 0}^4,$ we  denote 
$\ell_1! \ell_2! \ell_3! \ell_4!$ by $\ell !$, 
$\xi_1^{\ell_1}\xi_2^{\ell_2}\xi_3^{\ell_3}\xi_4^{\ell_4}$ by $\xi^\ell$,   
$\partial_1^{\ell_1} \partial_2^{\ell_2}
\partial_3^{\ell_3}\partial_4^{\ell_4}$ by $\partial^\ell $,  $\delta_1^{\ell_1} 
\delta_2^{\ell_2} \delta_3^{\ell_3} \delta_4^{\ell_4}$ by $\delta^\ell$, and 
$U_1^{\ell_1}U_2^{\ell_2}U_3^{\ell_3}U_4^{\ell_4}$ by $U^\ell$.

\begin{definition} 

A smooth map $\rho: \mathbb{R}^4 \to C^\infty(\mathbb{T}_\theta^4)$ is  
said to be a symbol of order $m \in \mathbb{Z}$, if for any set of 
non-negative integers $i, j \in \mathbb{Z}_{\geq 0}^4$,  there 
exists a constant $c$ such that  
\[ 
|| 
\partial^j \delta^i \big( 
\rho(\xi) \big )|| \leq c (1+|\xi|)^{m-|j|},
\]
and if there exists a smooth map 
$k: \mathbb{R}^4 \setminus \{0 \}\to C^\infty(\mathbb{T}_\theta^4)$ such that
\[
\lim_{\lambda \to \infty} \lambda^{-m} \rho(\lambda \xi) = k (\xi), 
\qquad \xi \in \mathbb{R}^4 \setminus \{0 \}.
\]

\end{definition}

The space of symbols of order $m$ is denoted by $S_m$. 
The pseudodifferential operator $P_\rho : C^\infty(\mathbb{T}_\theta^4) \to 
C^\infty(\mathbb{T}_\theta^4)$ associated to a symbol $\rho \in S_m$ is 
given by 
\[ 
P_{\rho}(a) = (2 \pi)^{-4} \int \int e^{-is \cdot \xi} \rho(\xi) 
\alpha_s(a) \,ds \,d\xi, \qquad a \in C^\infty(\mathbb{T}_\theta^4), 
\] 
where $\{ \alpha_s\}_{s \in \mathbb{R}^4}$ is the group of 
$C^*$-algebra automorphisms defined by
\[ 
\alpha_s(U^\ell ) =
e^{is \cdot \ell} U^\ell, \qquad \ell \in \mathbb{Z}^4, 
\]
which was explained in more generality in Subsection \ref{ncrtori}. 
For example, the differential operator  
$\sum a_\ell \delta^{\ell} $ 
is associated with the symbol 
$ \sum a_\ell \xi^{\ell} $ 
via the above formula.

The pseudodifferential operators on $\mathbb{T}_\theta^4$ form 
an algebra and there is an asymptotic expansion for the symbol 
of the composition of two such operators. There is also an asymptotic 
formula for the symbol of the formal adjoint of a pseudodifferential 
operator, where the adjoint is taken with respect to the 
inner product given by \eqref{innerp-vol}. We explain this in 
the following proposition in which the relation 
$\rho \sim \sum_{j=0}^\infty \rho_j$ for given symbols $\rho, \rho_j$ 
means that for any $k \in \mathbb{Z}$, 
there exists an $N \in \mathbb{Z}_{\geq 0}$ such that 
$\rho - \sum_{j=0}^{n} \rho_j \in S_k$ for any $n >N$.

\begin{proposition} \label{symbolcalculus}

Let $\rho \in S_m$ and $\rho' \in S_{m'}$. There exists a 
unique $\lambda \in S_{m+m'}$ such that 
\[
P_\lambda = P_\rho P_{\rho'}.
\]
Moreover,  
\[ 
\lambda \sim \sum_{\ell \in \mathbb{Z}^4_{\geq 0} } 
\frac{1}{\ell ! }
\partial^{\ell} \rho (\xi) \,
\delta^{\ell} (\rho'(\xi)).
\]
There is also a unique $\tau \in S_m$ such that $P_\tau$ is the 
formal adjoint of $P_\rho$, and 
\[
\tau \sim \sum_{\ell \in \mathbb{Z}^4_{\geq 0}} 
\frac{1}{\ell !}
\partial^{\ell} \delta^{\ell} \big (\rho(\xi) \big)^*.
\]

\end{proposition}

Elliptic pseudodifferential operators on $\mathbb{T}_\theta^4$ are 
defined to be those whose symbols have the following property: 
\begin{definition}
A symbol $\rho \in S_m$ is said to be elliptic if $\rho(\xi)$ is 
invertible for any $\xi \neq 0$, and if 
there exists a constant $c$ such that
\[ || \rho(\xi)^{-1} || \leq c (1+|\xi|)^{-m}, \]
when $|\xi|$ is sufficiently large. 
\end{definition}
As an example, the flat Laplacian 
$\triangle=\delta_1^2+\delta_2^2+\delta_3^2+\delta_4^2$ is an 
elliptic operator of order $2$ since its symbol is 
$\rho(\xi)=\xi_1^2+\xi_2^2+\xi_3^2+\xi_4^2$, 
which satisfies the above criterion.

\subsection{Small time asymptotic expansion for $\textnormal{Trace}(e^{-t \triangle_\varphi})$.} \label{smtae}

Geometric invariants of a Riemannian manifold, such as its volume and scalar 
curvature, can be computed by considering small time heat kernel expansions 
of its Laplacian. Pseudodifferential calculus may be employed to compute the 
terms of such expansions (cf. \cite{gil}). One can use Connes' pseudodifferential 
calculus \cite{con1}  to apply similar heat kernel techniques in order to compute 
geometric invariants of noncommutative spaces \cite{contre, conmos2, fatkha1, 
fatkha2, fatkha3}. Here, we briefly explain this method and derive the small time 
asymptotic expansion for the trace of $e^{-t \triangle_\varphi}$, where 
$\triangle_\varphi$ is the perturbed Laplacian on $\mathbb{T}_\theta^4$  
introduced in Subsection \ref{perturbedlaplacian}. First we need to 
compute the symbol of this differential operator:

\begin{lemma}

The symbol of $\triangle_\varphi$ is equal to
$ a_2 (\xi)+ a_1(\xi) + a_0(\xi), $ where
\begin{eqnarray}
&& a_2 (\xi) =  e^h \sum_{i=1}^4  \xi_i^2, \qquad 
 a_1 (\xi)=  \sum_{i=1}^4  \delta_i(e^h) \xi_i,  \nonumber \\
&& a_0(\xi) =  \sum_{i=1}^4 \big (  \delta_i^2(e^h) - \delta_i(e^h)e^{-h} 
\delta_i(e^h) \big ). \nonumber
\end{eqnarray}

\begin{proof}

It follows easily from the symbol calculus  explained in 
Proposition \ref{symbolcalculus}.

\end{proof}
  
\end{lemma}

Using the Cauchy integral formula, one has
\begin{equation} \label{cauchyint}
e^{-t\triangle_\varphi} = \frac{1}{2\pi i} \int_C e^{-t \lambda} 
(\triangle_\varphi - \lambda)^{-1} \, d \lambda, 
\end{equation}
where $C$ is a curve in the complex plane that goes around the non-negative
real axis  in such a way that
\[ 
e^{-t s} = \frac{1}{2\pi i} \int_C e^{-t \lambda} (s - \lambda)^{-1} \, 
d \lambda, \qquad s \geq 0.
\]
Appealing to this formula, one can use Connes' pseudodifferential 
calculus to  employ similar arguments to those in \cite{gil} and derive an 
asymptotic expansion of the form 
\begin{equation} \label{asympl}
\text{Trace}(e^{-t \triangle_\varphi}) \sim t^{-2} \sum_{n=0}^{\infty} 
B_{2n} (\triangle_\varphi)
t^n \qquad (t \to 0). \nonumber
\end{equation} 
That is, one can approximate $(\triangle_\varphi - \lambda)^{-1}$ 
by pseudodifferential operators $B_\lambda$ whose symbols are of 
the form 
\[
b_0(\xi, \lambda) +  b_1(\xi, \lambda) +  b_2(\xi, \lambda) + \cdots, 
\]
where for $j=0, 1, 2, \dots$, $b_j(\xi, \lambda)$ is a symbol of order $-2-j$.

Therefore, we have to use the calculus of symbols explained in 
Subsection \ref{ncpseudos}  to solve the equation 
\begin{eqnarray} 
(b_0+b_1+b_2+\cdots) ((a_2-\lambda)+a_1+a_0) \sim 1. \nonumber
\end{eqnarray}
Here, $\lambda$ is treated as a symbol of order $2$ and we let 
$a_2'=a_2-\lambda, a_1'=a_1, a_0'=a_0$. Then the above equation yields
\[ 
\sum_{\substack{j \in \mathbb{Z}_{\geq 0}, \, \ell \in \mathbb{Z}^4_{\geq 0},\\ k=0, 1, 2}} 
\frac{1}{\ell! } \partial^{\ell} b_j \,
\delta^{\ell} (a_k') \sim 1.
\]
Comparing symbols of the same order on both sides, one concludes that 
\begin{equation} 
b_0=a_2'^{-1}=(a_2 - \lambda)^{-1}= \Big ( e^h \sum_{i=1}^4  
\xi_i^2 - \lambda \Big )^{-1}, \nonumber
\end{equation}
and
\begin{equation} \label{bnformula}
b_n  = - \sum \frac{1}
{\ell! } 
\partial^{\ell} b_j \,
\delta^{\ell} (a_k) b_0, 
\qquad n>0,
\end{equation}
where the summation is over all 
$ 0 \leq j <n,  0 \leq k \leq 2, \ell \in \mathbb{Z}^4_{\geq 0}$ 
such that $2+j+|\ell |-k=n.$ 
Similar to \cite{gil, contre} one can use these symbols to approximate 
$e^{-t \triangle_\varphi}$ with suitable infinitely 
smoothing operators and derive the desired asymptotic expansion. 
We record this result in the following proposition.

\begin{proposition} \label{asympexpa}

There is a small time asymptotic expansion 
\begin{equation} 
\textnormal{Trace}(e^{-t \triangle_\varphi}) \sim t^{-2} \sum_{n=0}^{\infty} 
B_{2n} (\triangle_\varphi) t^n \qquad (t \to 0), \nonumber
\end{equation}
where  for each $n=0, 1, 2, \dots,$ 
\[
B_{2n}(\triangle_\varphi) = \frac{1}{2 \pi i} \int \int_C e^{-\lambda} 
 \varphi_0(b_{2n}(\xi, \lambda)) \, d \lambda\, d \xi.
 \]

\end{proposition}

\section{Weyl's Law and Connes' Trace Theorem} \label{WLCTT}

A celebrated theorem of Weyl states that one can hear the volume 
of a closed Riemannian manifold $(M, g)$ from the asymptotic distribution 
of the eigenvalues of its Laplacian $\triangle_g$ acting on smooth functions on $M$. 
That is, if $0 \leq \lambda_0 \leq \lambda_1 \leq \lambda_2 \leq \cdots$ are the 
eigenvalues of $\triangle_g$ 
counted with multiplicity and $N(\lambda)=\#\{\lambda_j \leq \lambda\}$ is 
the eigenvalue counting function then 
\[ 
N(\lambda) \sim \frac{\textnormal{Vol}(M)}{(4 \pi)^{n/2}\Gamma(\frac{n}{2}+1) } 
\lambda^{n/2} \qquad (\lambda \to \infty),
\]
where $n= \textnormal{dim}\,M$ and $\textnormal{Vol}(M)$ is the 
volume of $M$. An equivalent formulation of this result is the following 
asymptotic estimate for the eigenvalues: 
\[ 
\lambda_j \sim \frac{4 \pi \Gamma(\frac{n}{2}+1)^{2/n}}
{\textnormal{Vol}(M)^{2/n}} j^{2/n} \qquad (j \to \infty). 
\]
This readily shows that $(1+\triangle_g)^{-n/2}$ is in the domain of the 
Dixmier trace  and 
\[
\textnormal{Tr}_\omega\big( (1+\triangle_g)^{-n/2} \big )= 
\frac{\textnormal{Vol}(M)}{4 \pi^{n/2} \Gamma(\frac{n}{2}+1)}.
\] 
A generalization of this result is Connes' trace theorem \cite{con2} 
which states that the Dixmier trace  and Wodzicki's noncommutative 
residue \cite{wod} coincide on  pseudodifferential operators of order $-n$ acting 
on the sections of a vector bundle over $M$. In the sequel we will provide more 
explanations about the Dixmier trace, the Wodzicki residue, and Connes' trace 
theorem.

In this section, we establish the analogue of Weyl's law for 
$\mathbb{T}_\theta^4$ by studying the asymptotic distribution of 
the eigenvalues of the Laplacian $\triangle_\varphi$. We will then prove 
the analogue of Connes' trace theorem for $\mathbb{T}_\theta^4$. 
This is done by introducing a noncommutative residue on the algebra 
of classical pseudodifferential operators on $\mathbb{T}_\theta^4$ and 
showing that it coincides with the Dixmier trace on the pseudodifferential 
operators of order $-4$. We closely follow the constructions and arguments 
given in \cite{fatkha3, fatwon} for the noncommutative 2-torus, and 
because of similarities in the arguments, we provide essentials of the 
proofs rather briefly.

\subsection{Asymptotic distribution of the eigenvalues of $\triangle_\varphi$.}

Let $0 \leq \lambda_0 \leq \lambda_1 \leq \lambda_2 \leq \cdots$ be  
the eigenvalues of $\triangle_\varphi$, counted with multiplicity. 
It follows from the asymptotic expansion 
\[
\textnormal{Trace}(e^{-t \triangle_\varphi}) = \sum_{j=0}^\infty e^{-t \lambda_j} 
\sim t^{-2} \sum_{n=0}^{\infty} B_{2n} (\triangle_\varphi) t^n \qquad (t \to 0),
\]
derived  in Proposition \ref{asympexpa}, that 
\[
\lim_{t\to 0^+} t^2 \sum e^{-t \lambda_j} = B_{0} (\triangle_\varphi). 
\] 
It readily follows from Karamata's Tauberian theorem \cite{bgv} 
that the corresponding eigenvalue counting function 
$N$ has the following asymptotic 
behavior: 
\[
N(\lambda) \sim \frac{B_{0} (\triangle_\varphi)}{\Gamma(3)} \lambda^2 
\qquad (\lambda \to \infty). 
\]
We establish the analogue of Weyl's law for $\mathbb{T}_\theta^4$ 
in the following theorem by computing $B_{0}(\triangle_\varphi)$ 
and observing that, up to a universal constant, it is equal to 
$\varphi(1)=\varphi_0(e^{-2h})$.

\begin{theorem} \label{weyllaw}

The eigenvalue counting function $N$ of the Laplacian 
$\triangle_\varphi$ on $\mathbb{T}_\theta^4$ satisfies

\begin{equation} \label{asympcount}
N(\lambda) \sim  \frac{ \pi^2 \varphi_0(e^{-2h})}{2}  
\lambda^2 \qquad (\lambda \to \infty).
\end{equation}

\begin{proof}

Using Proposition \ref{asympexpa}, we have 
\begin{eqnarray} 
B_0(\triangle_\varphi) &=& 
\frac{1}{2 \pi i} \int \int_C e^{-\lambda} 
\varphi_0 (b_0(\xi, \lambda)) \, d \lambda \, d \xi  \nonumber \\
&=& \varphi_0 \Big ( \frac{1}{2 \pi i} \int \int_C e^{-\lambda} 
\big ( e^h \sum_{i=1}^4  \xi_i^2 - \lambda \big )^{-1} \, 
d \lambda \, d \xi  \Big  ) \nonumber \\
&=& \varphi_0 \Big ( \int e^{- e^h( \xi_1^2 + 
\cdots + \xi_4^2)} \, d\xi\Big ) \nonumber \\ 
&=&  \pi^2 \varphi_0(e^{-2h}). \nonumber
\end{eqnarray}
Thus, it follows from the above discussion that 
\[
N(\lambda) \sim \frac{ \pi^2 \varphi_0(e^{-2h})}{\Gamma(3)} \lambda^2  
=\frac{\pi^2 \varphi_0(e^{-2h})}{2}  \lambda^2 \qquad (\lambda \to \infty).
\]

\end{proof}

\end{theorem}

A corollary of this theorem is that $(1+\triangle_\varphi)^{-2}$ is 
in the domain of the Dixmier trace. Before stating the corollary we 
quickly review the Dixmier trace and the  noncommutative integral, 
following  \cite{con3}.

We denote the ideal of compact operators on a  Hilbert space 
$\mathcal{H}$ by $\mathcal{K}(\mathcal{H})$. For  any 
$T \in \mathcal{K}(\mathcal{H})$,  let $\mu_n(T), n=1, 2, \dots,$ 
denote the sequence of eigenvalues of its absolute value 
$|T|=(T^*T)^{\frac{1}{2}}$ written in  decreasing order with multiplicity:
\[
\mu_1(T) \geq \mu_2(T) \geq \cdots  \geq 0.
\]
The Dixmier trace is a trace functional on an ideal of 
compact operators $\mathcal{L}^{1,\infty } (\mathcal{H})$ defined as
\[
\mathcal{L}^{1,\infty } (\mathcal{H})= \Big \{T \in \mathcal{K}(\mathcal{H}); 
\qquad \sum_{n=1}^N \mu_n (T)=O \,(\text{log} N ) \Big \}.
\]
This ideal of operators is equipped with a natural norm:
\[
||T||_{1, \infty} := \sup_{N \geq 2} \frac{1}{\log N}\sum_{n=1}^N 
\mu_n(T), \qquad T \in \mathcal{L}^{1,\infty } (\mathcal{H}).
\]
It is clear that trace class operators are automatically in 
$\mathcal{L}^{1,\infty} (\mathcal{H})$. The Dixmier trace of an 
operator $T \in \mathcal{L}^{1,\infty}(\mathcal{H})$ measures the 
logarithmic divergence of
its ordinary trace. More precisely,  for any positive operator 
$T\in \mathcal{L}^{1,\infty } (\mathcal{H})$ we are interested 
in  the limiting behavior  of the  sequence
\[
\frac{1}{\log N}\sum_{n=1}^N \mu_n(T), \qquad N = 2, 3, \dots \,.
\]
While by our assumption this sequence is bounded,  its usual 
limit may not exist and must be replaced by a suitable 
generalized limit. The limiting procedure is carried out by 
means of a state on a $C^*$-algebra. Recall that a state 
on a  $C^*$-algebra is a non-zero positive linear functional 
on the algebra.

To define the Dixmier trace of a positive operator 
$T \in \mathcal{L}^{1,\infty}(\mathcal{H})$, consider 
the partial trace
\[
\text{Trace}_{N} (T) =\sum_{n=1}^N \mu_n (T), \qquad N=1, 2, \dots,
\]
and its piecewise affine interpolation denoted by 
$\text{Trace}_{r} (T)$ for  $r \in [1, \infty).$ 
Then let 
\[
\tau_{\Lambda} (T): =  \frac{1}{\log  \Lambda} \int_e^{\Lambda} 
\frac{ \text{Trace}_{r} (T)}{\log r}\frac{dr}{r}, \qquad 
\Lambda \in [e, \infty), 
\]
be the  Ces\`{a}ro mean of the function $\text{Trace}_{r} (T)/\log r$ 
over the multiplicative group $\mathbb{R}_{>0}$.
Now choosing a normalized state $\omega : C_b[e, \infty) \to 
\mathbb{C}$  on the algebra of bounded continuous functions
on  $[e, \infty)$ such that $\omega (f)=0$ for all $f$ vanishing at 
$\infty$, the Dixmier trace of  $T \geq 0$ is defined as 
\[ 
\text{Tr}_{\omega}(T) = \omega (\tau_{\Lambda} (T)).
\]
Then one can extend $\text{Tr}_{\omega}$ to all of 
$\mathcal{L}^{1,\infty}(\mathcal{H})$ by linearity. 

The resulting linear functional  $\text{Tr}_{\omega}$ is a 
positive trace on $\mathcal{L}^{(1, \infty)} (\mathcal{H})$ which 
in general depends on the limiting procedure $\omega$. 
The operators $T \in \mathcal{L}^{1,\infty}(\mathcal{H})$ 
whose Dixmier trace $\text{Tr}_{\omega}(T)$ is independent 
of the choice of the state $\omega$ are called 
measurable and we will denote their Dixmier trace by 
$ \int\!\!\!\!\!\!- \, T$. If for a compact positive operator $T$ we have 
\[
\mu_n(T) \sim \frac{c}{n}  \qquad (n \to \infty), 
\] 
where $c$ is a constant, then $T$ is measurable and $\int\!\!\!\!\!\!- \,T = c.$
We use this fact in the proof of the following corollary of Theorem \ref{weyllaw}.

\begin{corollary} \label{indomdix}

The operator $(1+\triangle_\varphi)^{-2}$, where $\triangle_\varphi$ is the 
perturbed Laplacian on $\mathbb{T}_\theta^4$, 
is a measurable operator in $\mathcal{L}^{1, \infty}(\mathcal{H}_0)$, 
and 
\[
\int\!\!\!\!\!\!- \,(1+\triangle_\varphi)^{-2}  =  \frac{\pi^2}{2} \varphi_0(e^{-2h}).
\]

\begin{proof}

It follows from the asymptotic behavior \eqref{asympcount}  of the 
eigenvalue counting function $N$ that the eigenvalues $\lambda_j$ of 
$\triangle_\varphi$ satisfy 
\[
\lambda_j \sim \frac{\sqrt{2}}{\pi \varphi_0(e^{-2h})^{1/2}} j^{1/2} 
\qquad (j \to \infty).
\]  
Therefore, using the above fact, it easily follows that 
$(1+\triangle_\varphi)^{-2}$ is measurable and 
$\int\!\!\!\!\!\!- \,(1+\triangle_\varphi)^{-2}  =  \frac{\pi^2}{2} \varphi_0(e^{-2h}).$

\end{proof}

\end{corollary}

\subsection{A noncommutative residue for $\mathbb{T}_\theta^4$.}

Let $M$ be a closed smooth manifold of dimension $n$. 
Wodzicki defined a trace functional on the algebra of 
pseudodifferential operators of arbitrary order on $M$, 
and proved that it was the only non-trivial trace \cite{wod}. 
This functional, denoted by $\text{Res}$, is called the 
noncommutative residue.

The noncommutative residue of a classical pseudodifferential 
operator $P$ acting on smooth sections of a vector bundle 
$E$ over $M$ is defined as 
\[
\text{Res}(P) = (2 \pi)^{-n} \int_{S^*M} 
{\rm tr}(\rho_{-n}(x, \xi))\, dx \, d\xi,
\]
where $S^*M \subset T^*M$ is the unit cosphere bundle 
on $M$ and $\rho_{-n}$ is the component of order $-n$ of 
the complete symbol of $P$.

Similar to \cite{fatwon},  we define a noncommutative residue 
on the algebra of classical pseudodifferential operators on 
$\mathbb{T}_\theta^4$, which is a natural analogue of the 
Wodzicki residue.

\begin{definition}

A pseudodifferential symbol $\rho \in S_m$ on $\mathbb{T}_\theta^4$ 
is said to be classical if there is an asymptotic expansion of the 
form 
\[
\rho(\xi) \sim \sum_{j=0}^{\infty} \rho_{m-j}(\xi) \qquad  (\xi \to \infty), 
\] 
where  each $\rho_{m-j} : \mathbb{R}^4 \setminus \{0 \} \to 
C^\infty(\mathbb{T}_\theta^4)$ is smooth and positively 
homogeneous of order $m-j$. Given such a symbol, we define the 
noncommutative residue of the corresponding pseudodifferential 
operator $P_\rho$ as 
\[
\textnormal{res}(P_\rho) = \int_{\mathbb{S}^3} 
\varphi_0 (\rho_{-4}(\xi)) \, d\Omega, 
\]
where $d\Omega$ is the standard invariant measure on the unit 
sphere in $\mathbb{R}^4$. 
\end{definition}

We note that, as shown in \cite{fatwon}, the homogeneous 
terms in the expansion of any classical pseudodifferential 
symbol are uniquely determined. Thus, there is no ambiguity 
in the above definition.

In the following theorem we identify all continuous trace functionals 
on the algebra of classical pseudodifferential operators on 
$\mathbb{T}_\theta^4$. A linear functional on this algebra 
is said to be continuous if it vanishes on the operators 
whose symbols are of sufficiently small order. 
First we state two lemmas which will be used in the proof 
of the following theorem. One can prove these lemmas 
by similar arguments to those given in \cite{fatwon}.

\begin{lemma} \label{homogder}

Let $f: \mathbb{R}^4 \setminus \{0\} \to C^\infty(\mathbb{T}_\theta^4)$ 
be a smooth map which is positively homogeneous
of order $m \in \mathbb{Z}$. If $m \neq -4 $, or if 
$m =-4$ and   $\int_{\mathbb{S}^3} f \,d\Omega =0,$ 
then one can write 
$
f = \sum_{i=1}^4\partial_i (h_i), 
$
for some smooth maps 
$h_i:\mathbb{R}^4 \setminus \{0\} \to C^\infty(\mathbb{T}_\theta^4)$.  

\end{lemma}

\begin{lemma} \label{sumsymbols}

Let $\sigma_j \in S_{m_j}$, $j=0, 1, 2, \dots$, be a 
sequence of pseudodifferential symbols on 
$\mathbb{T}_\theta^4$ such that $\lim m_j = - \infty$. 
There exists a symbol $\sigma$ such that
$ \sigma \sim \sum_{j=0}^{\infty} \sigma_{j}.$

\end{lemma}

\begin{theorem}

The noncommutative residue $\textnormal{res}$ is a trace, and 
up to multiplication by a constant, it is the unique continuous 
trace on the algebra of classical pseudodifferential operators 
on $\mathbb{T}_\theta^4$. 

\begin{proof}

Let $\rho, \rho': \mathbb{R}^4 \to C^\infty(\mathbb{T}_\theta^4) $ 
be classical symbols with asymptotic expansions
\[ 
\rho(\xi) \sim \sum_{j=0}^{\infty} \rho_{n-j}(\xi), \qquad  
\rho'(\xi) \sim \sum_{k=0}^{\infty} \rho_{n'-k}(\xi) \qquad  (\xi \to \infty),
\]
where  $\rho_{n-j}$ is homogeneous of order $n-j$ and $\rho_{n'-k}$ is 
homogeneous of order $n'-k$. Using the calculus of symbols explained 
in Proposition \ref{symbolcalculus} and the trace property of $\varphi_0$, 
we have:
\begin{eqnarray}
 &&\textnormal{res}(P_\rho P_{\rho'} - P_{\rho'} P_\rho) \nonumber \\
&&=\int_{\mathbb{S}^3}  \varphi_0 \Big (            
\sum \frac{1}{\ell !} \big ( 
\partial^\ell (\rho_{n-j}) \delta^\ell (\rho'_{n'-k})
- \delta^\ell (\rho_{n-j}) \partial^\ell (\rho'_{n'-k})
\big ) 
\Big ) \, d \Omega, \nonumber
\end{eqnarray}
where the summation is over all $j, k \in \mathbb{Z}_{\geq 0}$ and 
$\ell \in  \mathbb{Z}^4_{\geq 0}$ such that $n+n'-j-k-|\ell |=-4$.  
One can write each $\partial^\ell (\rho_{n-j}) \delta^\ell (\rho'_{n'-k})
- \delta^\ell (\rho_{n-j}) \partial^\ell (\rho'_{n'-k})$ in the 
above integral as 
\[
\sum_{i=1}^4 \big ( \partial_i(f_i) + \delta_i(g_i) \big ), 
\]
for some $f_i, g_i: \mathbb{R}^4 \setminus \{ 0 \} \to 
C^\infty(\mathbb{T}_\theta^4)$, where each $f_i$ is 
positively homogeneous of order $-3$. Thus, using 
Lemma  5.1.3 on page 208 of \cite{nicrod} and 
the fact that $\varphi_0 \circ \delta_i =0$, 
we have 
\[
\textnormal{res}(P_\rho P_{\rho'} - P_{\rho'} P_\rho)=0.
\] 
This proves the trace property of $\textnormal{res}$. 

In order to prove the uniqueness, assume that $\psi$ 
is a continuous trace on the algebra of classical symbols 
on $\mathbb{T}_\theta^4$. For any classical symbol 
$\rho$,  the symbol of $P_{\xi_i} P_\rho - P_\rho P_{\xi_i}$ 
is equivalent to $\delta_i(\rho)$. Since $\psi$ is a trace, it 
follows that 
\begin{equation} \label{trdeltazero}
\psi(P_{\delta_i(\rho)})=0.
\end{equation}

Since we can write 
\[
\rho=\varphi_0(\rho)+\sum_{i=1}^4\delta_i(\rho_i),
\]
for some symbols $\rho_i$, equation \eqref{trdeltazero} 
implies that 
\begin{equation} \label{trconstpart}
\psi(P_\rho)=\psi(P_{\varphi_0(\rho)}).
\end{equation}

Also, by considering the symbol of 
$P_\rho P_{U_i} - P_{U_i} P_\rho$ 
and using similar arguments to those in \cite{fatwon}, one 
can conclude that 
\begin{equation} \label{trzeropartial}
\psi(P_{\partial_i(\rho)})=0.
\end{equation}

Now we consider the asymptotic expansion 
\[
\rho(\xi) \sim \sum_{j=0}^\infty \rho_{n-j}(\xi), \qquad (\xi \to \infty), 
\]
where $\rho_{n-j}$ is positively homogeneous of order $n-j$. 
By using Lemmas \ref{homogder} and \ref{sumsymbols}, and 
setting 
$r=\textnormal{Vol}(\mathbb{S}^3)^{-1} \int_{\mathbb{S}^3} 
\rho_{-4} \, d\Omega ,$ we can write 
\begin{eqnarray} \label{symbolexpdd}
\rho &\sim& \rho_{-4} + \sum_{n-j \neq -4} \rho_{n-j} \nonumber \\
&=& \frac{r}{|\xi|^4}+ \big (\rho_{-4}- \frac{r}{|\xi|^4} \big )+ 
\sum_{n-j \neq -4} \sum_{i=1}^4 \partial_i (\rho_{n-j,i}) \nonumber \\
&\sim& \frac{r}{|\xi|^4}+\sum_{i=1}^4 \partial_i (\rho_{-4,i}) +
\sum_{i=1}^4 \partial_i \Big ( \sum_{n-j \neq -4} \rho_{n-j,i} \Big ), 
\end{eqnarray}
for some smooth maps 
$\rho_{-4, i}, \rho_{n-j,i}: \mathbb{R}^4 \setminus \{ 0\} \to 
C^\infty(\mathbb{T}_\theta^4)$. Now we can use \eqref{trconstpart}, 
\eqref{trzeropartial}, \eqref{symbolexpdd} 
to conclude that 
\[
\psi(P_\rho) = \psi(P_{r/|\xi^4|})=\psi(P_{\varphi_0(r/|\xi^4|)})=
\varphi_0(r) \psi(P_{1/|\xi^4|})=
\frac{\psi(P_{1/|\xi^4|})}{\textnormal{Vol}(\mathbb{S}^3)} 
\textnormal{res}(P_\rho). 
\]

\end{proof}

\end{theorem}

\subsection{Connes' trace theorem for $\mathbb{T}_\theta^4$.}

As above, let $M$ be a closed smooth manifold of dimension 
$n$. The restriction of the Wodzicki residue $\text{Res}$ to 
pseudodifferential operators of order $-n$ was discovered 
independently by Guillemin and its properties were studied in 
\cite{guil}. In general, unlike the Dixmier trace, $\text{Res}$ 
is not a positive linear functional. However,  
its restriction to pseudodifferential operators of order $-n$ 
is positive.  One of the main results proved in \cite{con2} 
is that if $E$ is a smooth vector bundle on $M$  then the 
Dixmier trace $\text{Tr}_\omega$ and $\text{Res}$ coincide 
on pseudodifferential operators of order $-n$  acting on 
$L^2$ sections of  $E$. In fact it is proved that such operators $P$ 
are measurable operators in $\mathcal{L}^{1,\infty}(L^2(M, E))$ and 
\[
 \int\!\!\!\!\!\!- \, P = \frac{1}{n} \text{Res}(P). 
\]
This result is known as Connes' trace theorem. In the following 
theorem, we establish the analogue of this result for the 
noncommutative $4$-torus $\mathbb{T}_\theta^4$.

\begin{theorem}

Let $\rho$ be a classical pseudodifferential symbol of order 
$-4$ on  $\mathbb{T}_\theta^4$. Then 
$P_\rho$ is a measurable operator in  $\mathcal{L}^{1, \infty}(\mathcal{H}_0),$ 
and under the assumption that all nonzero entries of $\theta$ are irrational, 
we have
\[
\int\!\!\!\!\!\!- \, P_\rho = 
\frac{1}{4} \, \textnormal{res}(P_\rho). 
\]

\begin{proof}

In order to show that $P_\rho \in \mathcal{L}^{1, \infty}(\mathcal{H}_0)$, 
we write $P_\rho=A(1+\triangle)^{-2}$, where 
$\triangle=\delta_1^2+\cdots+\delta_4^2$ is the flat Laplacian, and 
$A=P_\rho(1+\triangle)^{2}$ is a pseudodifferential operator of 
order 0. Since $A$ is a bounded operator on $\mathcal{H}_0$ and it was 
shown in Corollary \ref{indomdix} that 
$(1+\triangle)^{-2} \in \mathcal{L}^{1, \infty}(\mathcal{H}_0)$, it follows 
that $P_\rho$ is in the domain of the Dixmier trace. Using a similar argument, 
one can see that  any pseudodifferential operator of order $-5$ on 
$\mathbb{T}_\theta^4$ is in the kernel of the Dixmier trace. Therefore, 
if we write 
\[ 
\rho(\xi) \sim \rho'(\xi) + \sum_{j=1}^{\infty} \rho_{-4-j}(\xi)  
\qquad  (\xi \to \infty),
\]
where $\rho'$ and $\rho_{-4-j}$ are respectively positively homogeneous 
of order $-4$ and $-4-j$,  then 
\begin{equation} \label{only-4}
{\rm Tr}_\omega (P_\rho) = {\rm Tr}_\omega (P_{\rho'}).
\end{equation}

Also, since the symbol of $P_{\rho'} P_{U_i} - P_{U_i} P_{\rho'}$ is equivalent 
to the symbol of $P_{\rho'U_i} -P_{U_i \rho'}$ modulo a symbol of order $-5$, 
and ${\rm Tr}_\omega$ is a trace, for $i=1, \dots 4,$ we have 
\[
{\rm Tr}_\omega(P_{\rho'U_i})={\rm Tr}_\omega(P_{U_i\rho'}).   
\]
It follows from this observation that if $f: \mathbb{R}^4\setminus \{ 0 \}
\to \mathbb{C}$ is smooth and positively homogeneous of order -4, then 
${\rm Tr}_\omega(P_{fU^\ell})=0$ if $\ell \neq 0 \in \mathbb{Z}^4$. 
Therefore, similar to the argument given in \cite{fatkha3} one can write the 
following expansion with rapidly decreasing coefficients 
\[
\rho'(\xi)=\sum_{\ell \in \mathbb{Z}^4} \rho'_\ell(\xi) U^\ell, 
\]
and conclude that 
\begin{equation} \label{onlytr-4}
{\rm Tr}_\omega(P_{\rho'})= {\rm Tr}_\omega (P_{\rho'_0}), 
\end{equation}
where, as noted above, $\rho'_0=\varphi_0 \circ \rho': 
\mathbb{R}^4\setminus \{ 0 \} \to \mathbb{C}.$ Note that one has to use 
the fact that ${\rm Tr}_\omega$ is continuous with respect to the uniform 
norm of symbols of order $-4$, namely that, if we let 
$q_m(\xi) = \rho'(\xi)- \sum_{|\ell|\leq m } a_\ell(\xi) U^\ell,$ then 
\[
\lim_{m\to \infty} {\rm Tr}_\omega (P_{q_m})=0.
\]

Since $\rho'_0$ is a smooth complex-valued function on 
$\mathbb{R}^4 \setminus \{0 \}$ which is positively homogeneous 
of order $-4$, in order to analyze ${\rm Tr}_\omega (P_{\rho'_0})$, we 
define a linear functional $\mu$ on the space of continuous complex-valued 
functions on $\mathbb{S}^3$, as follows. Given a smooth function 
$f:\mathbb{S}^3 \to \mathbb{C}$, we denote its positively homogeneous 
extension of order $-4$ to $\mathbb{R}^4 \setminus \{0 \}$  by $\tilde{f}$, and 
define $\mu(f)={\rm Tr}_\omega(P_{\tilde{f}})$. Using the continuity property 
mentioned above, $\mu$ extends to the space of continuous functions on 
$\mathbb{S}^3$. Also using positivity of ${\rm Tr}_\omega$, one can see 
that $\mu$ is a positive linear functional. Thus, it follows from the Reisz 
representation theorem that $\mu$ is given by integration against a Borel 
measure on $\mathbb{S}^3$. This measure is rotation invariant, which 
can be shown by  using the trace property of ${\rm Tr}_\omega$ and the 
fact that for any rotation $T$ of $\mathbb{R}^4$ and any pseudodifferential 
symbol $\sigma$ on $\mathbb{R}^4$ we have 
\[
P_{\sigma(Tx, T\xi)}=\mathcal{U}^{-1}P_{\sigma(x, \xi)} \mathcal{U}, 
\] 
where $\mathcal{U}$ is the unitary operator $\mathcal{U}(g)=g \circ T^{-1}$, 
$g \in C_c^\infty(\mathbb{R}^4)$. Therefore $\mu$ is given by integration 
against a constant multiple of the standard invariant measure on $\mathbb{S}^3$. 
Denoting this constant by $c$, identities \eqref{only-4}, \eqref{onlytr-4} imply 
that 
\[
{\rm Tr}_\omega(P_\rho)={\rm Tr}_\omega(P_{\rho'})=
{\rm Tr}_\omega(P_{\rho'_0})= \mu(\rho'_0{\big \arrowvert}_{\mathbb{S}^3})
= c \int_{\mathbb{S}^3} \rho'_0{\big \arrowvert}_{\mathbb{S}^3} \, d\Omega 
= c \,{\rm res}(P_\rho).  
\]

The constant $c$ can be fixed by considering the flat Laplacian 
$\triangle=\sum_{i=1}^4 \delta_i^2$. According to corollary \ref{indomdix}, we have 
${\rm Tr}_\omega((1+\triangle)^{-2})= \pi^2/2$. On the other hand, considering 
the fact that the term of order $-4$ in the asymptotic expansion of the symbol of 
$(1+\triangle)^{-2}$ is $|\xi|^{-4}$, we have ${\rm res}((1+\triangle)^{-2})= 2 \pi^2$. 
Therefore $c=1/4$. 

\end{proof}

\end{theorem}

\section{Scalar Curvature and Einstein-Hilbert Action} \label{SCEHA}

Let $(M, g)$ be a smooth compact manifold of dimension $n \geq 2,$ 
and $\triangle_g$ be the Laplacian acting on smooth functions on $M$. 
For any $t >0$, the operator $e^{-t \triangle_g}$ is an infinitely smoothing 
operator, and there is an asymptotic expansion for its kernel $K(t, x, y)$, 
which is of the form 
\[
K(t, x, y) \sim  \frac{e^{- \textrm{dist}(x, y)^2/4t}}{(4 \pi t)^{n/2}}
\big (u_0(x, y)+ u_1(x, y)t + u_2(x, y)t^2+\cdots \big )  
\qquad (t \to 0).   
\] 
The coefficients $u_i$ are smooth functions defined in a neighborhood of 
the diagonal in $M\times M$. The kernel $K$ is called the heat kernel 
since it is the fundamental solution of the heat operator $\partial_t + \triangle_g$, 
and the $u_i$ are called the heat kernel coefficients. These coefficients are 
well-studied and are often known under the names of people who made the major 
contributions in the study, namely that, they are called both 
Minakshisundaram-Pleijel coefficients and Seeley-De Witt coefficients. 
Minakshisundaram and Pleijel derived the above asymptotic expansion 
by using the transport equation method of Hadamard. An approach, which is 
followed in the noncommutative case in \cite{contre, conmos2, fatkha1, fatkha2, fatkha3}, 
is to use pseudodifferential calculus to derive such asymptotic expansions. For a 
clear account of this approach and a detailed discussion of the local geometric 
information that are encoded in heat coefficients, we refer the reader to \cite{gil} 
and the references therein.

A crucial point that is used to define and compute the scalar curvature for 
noncommutative spaces \cite{conmar, conmos2, fatkha2} is that, up to a 
universal factor, the restriction of $u_1$ to the diagonal gives the scalar 
curvature of $M$. This, via the Mellin transform, allows to have a spectral 
definition for scalar curvature in terms of values or residues of spectral 
zeta functions. That is, if for instance the dimension of $M$ is $4$, the scalar 
curvature is the unique  $R \in C^\infty(M)$ (up to a universal constant) such 
that 
\[
\textrm{res}_{s=1} \textrm{Trace} (f \triangle_g^{-s}) = \int_M f R \,dvol_g, 
\]
for any $f \in C^\infty(M)$.

In this section we define the scalar curvature for $\mathbb{T}_\theta^4$ 
equipped with the perturbed Laplacian $\triangle_\varphi$, and compute 
the functions that give a local expression for the curvature. Then we 
consider the analogue of the Einstein-Hilbert action $\int_M R \, dvol_g$. 
We find a local expression for this action as well, and  show 
that its extremum is attained if and only if the Weyl factor is a constant, 
which is equivalent to having a metric with constant curvature.

\subsection{Scalar curvature for $\mathbb{T}_\theta^4$.} \label{SCsub}

Following \cite{conmar, conmos2, fatkha2} and the above discussion, we 
define the scalar curvature  of $\mathbb{T}_\theta^4$ equipped with the 
perturbed Laplacian $\triangle_\varphi$ as follows 
(see also \cite{bhumar} for a variant). 

\begin{definition}

The scalar curvature of the noncommutative 4-torus equipped with 
the perturbed volume form is the unique element 
$R \in C^\infty(\mathbb{T}_\theta^4)$ such that 
\begin{equation} 
\textnormal{res}_{s=1} \textnormal{Trace}(a \triangle_{\varphi}^{-s}) = 
\varphi_0(a R), \nonumber
\end{equation}
for any $a \in C^\infty(\mathbb{T}_\theta^4)$. 

\end{definition}

We follow the method employed in \cite{contre, conmos2, fatkha2} to 
find a local expression for the scalar curvature $R$. For the sake of 
completeness, we explain this procedure in the following proposition.

\begin{proposition}
The scalar curvature $R$ is equal to 
\[  
\frac{1}{2 \pi i} \int_{\mathbb{R}^2} \int_C 
e^{-\lambda} b_2(\xi, \lambda) \, d\lambda \, d\xi ,    
\]
where $b_2$ is the term of order $-4$ of the pseudodifferential 
symbol of the parametrix of $\triangle_\varphi-\lambda$, given in 
Subsection \ref{smtae}.
\begin{proof}
Using the Mellin transform we have
\[
a \triangle_\varphi^{-s} = \frac{1}{\Gamma(s)} \int_0^\infty 
a (e^{-t \triangle_\varphi} - P)t^{s-1} \,dt, \qquad 
a \in C^\infty(\mathbb{T}_\theta^4), \qquad \Re(s) >0, 
\]
where $P$ denotes the orthogonal projection on $\text{Ker}(\triangle_\varphi)$.

Appealing to the Cauchy integral formula \eqref{cauchyint} and 
using similar arguments to those explained in Subsection \ref{smtae} 
(cf. \cite{gil, contre}), one can derive an asymptotic expansion of the form 
\[ 
\text{Trace}(ae^{-t \triangle_\varphi}) \sim 
t^{-2} \sum_{n=0}^{\infty} B_{2n} (a, \triangle_\varphi)  t^n \qquad (t \to 0).  
\]
Using this asymptotic expansion  one can see that the zeta function 
\[
\zeta_a(s) = \textnormal{Trace}(a \triangle_\varphi^{-s}), 
\qquad \Re(s) \gg 0, 
\]
has a meromorphic extension to the whole plane with a simple pole at $1$, 
and
\[  
\textnormal{res}_{s=1}\zeta_a(s) = B_2(a, \triangle_\varphi).
\]
On the other hand, there are explicit formulas for the coefficients 
of the above asymptotic expansion. In particular we have
\begin{eqnarray}
B_2(a, \triangle_\varphi) &=& \frac{1}{2 \pi i} \int \int_C 
e^{-\lambda}  \varphi_0 \big (a b_2(\xi, \lambda) \big )\, d \lambda 
\, d \xi. \nonumber 
\end{eqnarray}

\end{proof}

\end{proposition}

We directly compute $b_2$ using \eqref{bnformula} and in order 
to  compute
\[
\frac{1}{ 2 \pi i} \int_{\mathbb{R}^4} \int_{C} e^{-\lambda} 
b_2(\xi, \lambda) \, d\lambda \, d\xi, 
\]
we use a homogeneity argument for the contour integral 
(cf. \cite{contre, conmos2})   and pass to the spherical 
coordinates   
\begin{eqnarray}
\xi_1 = r \sin (\psi) \sin (\phi) \cos (\theta), && \qquad   
\xi_2 =r \sin (\psi) \sin (\phi) \sin(\theta),  \nonumber \\
\xi_3 =r \sin (\psi) \cos (\phi), \,\, \,\, \quad  \quad       &&       
\qquad \xi_4 = r \cos(\psi),      \nonumber
\end{eqnarray}  
with $0 \leq r < \infty,$ $0 \leq \psi < \pi,$ $0 \leq \phi < \pi,$ 
$0 \leq \theta < 2 \pi.$ After working out 
the integrations with respect to the angles $\psi, \phi, \theta$,  
we obtain the following terms up to 
an overall factor of $-\pi^2$:  \\

\noindent 
$
+4r^9e^h b_0 b_0 b_0 \delta_1(e^h) e^h b_0 b_0 \delta_1(e^h) b_0 
+2r^9 e^h b_0 b_0 e^h b_0 b_0 \delta_1(e^h) b_0 \delta_1(e^h) b_0  \\
+4r^9 e^h b_0 b_0 \delta_1(e^h) e^h b_0 b_0 b_0 \delta_1(e^h) b_0  
+2r^9 e^h b_0 b_0 \delta_1(e^h) e^h b_0 b_0 \delta_1(e^h) b_0 b_0  \\
+4r^9 e^h b_0 e^h b_0 b_0 b_0 \delta_1(e^h) b_0 \delta_1(e^h) b_0  
+2r^9 e^h b_0 e^h b_0 b_0 \delta_1(e^h) b_0 b_0 \delta_1(e^h) b_0  \\
+2r^9 e^h b_0 e^h b_0 b_0 \delta_1(e^h) b_0 \delta_1(e^h) b_0 b_0  
+6r^9 e^h e^h b_0 b_0 b_0 b_0 \delta_1(e^h) b_0 \delta_1(e^h) b_0  \\
+2r^9 e^h e^h b_0 b_0 b_0 \delta_1(e^h) b_0 b_0 \delta_1(e^h) b_0  
+2r^9 e^h e^h b_0 b_0 b_0 \delta_1(e^h) b_0 \delta_1(e^h) b_0 b_0  \\
-r^7 b_0 b_0 \delta_1(e^h)e^h b_0 b_0 \delta_1(e^h) b_0 
-2r^7 b_0 \delta_1(e^h) e^h b_0 b_0 b_0 \delta_1(e^h) b_0  \\
-r^7 b_0 \delta_1(e^h) e^h b_0 b_0 \delta_1(e^h) b_0 b_0  
-16r^7 e^h b_0 b_0 b_0 \delta_1(e^h) b_0 \delta_1(e^h) b_0  \\
-r^7 e^h b_0 b_0 e^h b_0 b_0 \delta_1 \delta_1(e^h) b_0  
-8r^7 e^h b_0 b_0 \delta_1(e^h) b_0 b_0 \delta_1(e^h) b_0  \\
-8r^7 e^h b_0 b_0 \delta_1(e^h) b_0 \delta_1(e^h) b_0 b_0  
-2r^7 e^h b_0 e^h b_0 b_0 b_0 \delta_1 \delta_1(e^h) b_0  \\
-r^7 e^h b_0 e^h b_0 b_0 \delta_1 \delta_1(e^h) b_0 b_0  
-3r^7 e^h e^h b_0 b_0 b_0 b_0 \delta_1 \delta_1(e^h) b_0  \\
-r^7 e^h e^h b_0 b_0 b_0 \delta_1 \delta_1(e^h) b_0 b_0  
+5/2r^5 b_0b_0 \delta_1(e^h) b_0 \delta_1(e^h) b_0  \\
+2r^3 b_0 b_0 \delta_1(e^h) e^{-h} \delta_1(e^h) b_0  
+5/2r^5 b_0 \delta_1(e^h) b_0 b_0 \delta_1(e^h) b_0 \\ 
+5/2r^5 b_0 \delta_1(e^h) b_0 \delta_1(e^h) b_0 b_0  
+2r^3 b_0 \delta_1(e^h) e^{-h} \delta_1(e^h) b_0 b_0 \\ 
+6r^5 e^h b_0 b_0 b_0 \delta_1 \delta_1(e^h) b_0  
+3r^5 e^h b_0 b_0 \delta_1 \delta_1(e^h) b_0 b_0  \\
-2r^3 b_0 b_0 \delta_1 \delta_1(e^h) b_0 
-2r^3b_0 \delta_1 \delta_1(e^h) b_0 b_0\\
+4r^9e^h b_0 b_0 b_0 \delta_2(e^h) e^h b_0 b_0 \delta_2(e^h) b_0 
+2r^9 e^h b_0 b_0 e^h b_0 b_0 \delta_2(e^h) b_0 \delta_2(e^h) b_0  \\
+4r^9 e^h b_0 b_0 \delta_2(e^h) e^h b_0 b_0 b_0 \delta_2(e^h) b_0  
+2r^9 e^h b_0 b_0 \delta_2(e^h) e^h b_0 b_0 \delta_2(e^h) b_0 b_0  \\
+4r^9 e^h b_0 e^h b_0 b_0 b_0 \delta_2(e^h) b_0 \delta_2(e^h) b_0  
+2r^9 e^h b_0 e^h b_0 b_0 \delta_2(e^h) b_0 b_0 \delta_2(e^h) b_0  \\
+2r^9 e^h b_0 e^h b_0 b_0 \delta_2(e^h) b_0 \delta_2(e^h) b_0 b_0  
+6r^9 e^h e^h b_0 b_0 b_0 b_0 \delta_2(e^h) b_0 \delta_2(e^h) b_0  \\
+2r^9 e^h e^h b_0 b_0 b_0 \delta_2(e^h) b_0 b_0 \delta_2(e^h) b_0  
+2r^9 e^h e^h b_0 b_0 b_0 \delta_2(e^h) b_0 \delta_2(e^h) b_0 b_0  \\
-r^7 b_0 b_0 \delta_2(e^h)e^h b_0 b_0 \delta_2(e^h) b_0 
-2r^7 b_0 \delta_2(e^h) e^h b_0 b_0 b_0 \delta_2(e^h) b_0  \\
-r^7 b_0 \delta_2(e^h) e^h b_0 b_0 \delta_2(e^h) b_0 b_0  
-16r^7 e^h b_0 b_0 b_0 \delta_2(e^h) b_0 \delta_2(e^h) b_0  \\
-r^7 e^h b_0 b_0 e^h b_0 b_0 \delta_2 \delta_2(e^h) b_0  
-8r^7 e^h b_0 b_0 \delta_2(e^h) b_0 b_0 \delta_2(e^h) b_0  \\
-8r^7 e^h b_0 b_0 \delta_2(e^h) b_0 \delta_2(e^h) b_0 b_0  
-2r^7 e^h b_0 e^h b_0 b_0 b_0 \delta_2 \delta_2(e^h) b_0  \\
-r^7 e^h b_0 e^h b_0 b_0 \delta_2 \delta_2(e^h) b_0 b_0  
-3r^7 e^h e^h b_0 b_0 b_0 b_0 \delta_2 \delta_2(e^h) b_0  \\
-r^7 e^h e^h b_0 b_0 b_0 \delta_2 \delta_2(e^h) b_0 b_0  
+5/2r^5 b_0b_0 \delta_2(e^h) b_0 \delta_2(e^h) b_0  \\
+2r^3 b_0 b_0 \delta_2(e^h) e^{-h} \delta_2(e^h) b_0  
+5/2r^5 b_0 \delta_2(e^h) b_0 b_0 \delta_2(e^h) b_0 \\ 
+5/2r^5 b_0 \delta_2(e^h) b_0 \delta_2(e^h) b_0 b_0  
+2r^3 b_0 \delta_2(e^h) e^{-h} \delta_2(e^h) b_0 b_0 \\ 
+6r^5 e^h b_0 b_0 b_0 \delta_2 \delta_2(e^h) b_0  
+3r^5 e^h b_0 b_0 \delta_2 \delta_2(e^h) b_0 b_0  \\
-2r^3 b_0 b_0 \delta_2 \delta_2(e^h) b_0 
-2r^3b_0 \delta_2 \delta_2(e^h) b_0 b_0 \\
+4r^9e^h b_0 b_0 b_0 \delta_3(e^h) e^h b_0 b_0 \delta_3(e^h) b_0 
+2r^9 e^h b_0 b_0 e^h b_0 b_0 \delta_3(e^h) b_0 \delta_3(e^h) b_0  \\
+4r^9 e^h b_0 b_0 \delta_3(e^h) e^h b_0 b_0 b_0 \delta_3(e^h) b_0  
+2r^9 e^h b_0 b_0 \delta_3(e^h) e^h b_0 b_0 \delta_3(e^h) b_0 b_0  \\
+4r^9 e^h b_0 e^h b_0 b_0 b_0 \delta_3(e^h) b_0 \delta_3(e^h) b_0  
+2r^9 e^h b_0 e^h b_0 b_0 \delta_3(e^h) b_0 b_0 \delta_3(e^h) b_0  \\
+2r^9 e^h b_0 e^h b_0 b_0 \delta_3(e^h) b_0 \delta_3(e^h) b_0 b_0  
+6r^9 e^h e^h b_0 b_0 b_0 b_0 \delta_3(e^h) b_0 \delta_3(e^h) b_0  \\
+2r^9 e^h e^h b_0 b_0 b_0 \delta_3(e^h) b_0 b_0 \delta_3(e^h) b_0  
+2r^9 e^h e^h b_0 b_0 b_0 \delta_3(e^h) b_0 \delta_3(e^h) b_0 b_0  \\
-r^7 b_0 b_0 \delta_3(e^h)e^h b_0 b_0 \delta_3(e^h) b_0 
-2r^7 b_0 \delta_3(e^h) e^h b_0 b_0 b_0 \delta_3(e^h) b_0  \\
-r^7 b_0 \delta_3(e^h) e^h b_0 b_0 \delta_3(e^h) b_0 b_0  
-16r^7 e^h b_0 b_0 b_0 \delta_3(e^h) b_0 \delta_3(e^h) b_0  \\
-r^7 e^h b_0 b_0 e^h b_0 b_0 \delta_3 \delta_3(e^h) b_0  
-8r^7 e^h b_0 b_0 \delta_3(e^h) b_0 b_0 \delta_3(e^h) b_0  \\
-8r^7 e^h b_0 b_0 \delta_3(e^h) b_0 \delta_3(e^h) b_0 b_0  
-2r^7 e^h b_0 e^h b_0 b_0 b_0 \delta_3 \delta_3(e^h) b_0  \\
-r^7 e^h b_0 e^h b_0 b_0 \delta_3 \delta_3(e^h) b_0 b_0  
-3r^7 e^h e^h b_0 b_0 b_0 b_0 \delta_3 \delta_3(e^h) b_0  \\
-r^7 e^h e^h b_0 b_0 b_0 \delta_3 \delta_3(e^h) b_0 b_0  
+5/2r^5 b_0b_0 \delta_3(e^h) b_0 \delta_3(e^h) b_0  \\
+2r^3 b_0 b_0 \delta_3(e^h) e^{-h} \delta_3(e^h) b_0  
+5/2r^5 b_0 \delta_3(e^h) b_0 b_0 \delta_3(e^h) b_0 \\ 
+5/2r^5 b_0 \delta_3(e^h) b_0 \delta_3(e^h) b_0 b_0  
+2r^3 b_0 \delta_3(e^h) e^{-h} \delta_3(e^h) b_0 b_0 \\ 
+6r^5 e^h b_0 b_0 b_0 \delta_3 \delta_3(e^h) b_0  
+3r^5 e^h b_0 b_0 \delta_3 \delta_3(e^h) b_0 b_0  \\
-2r^3 b_0 b_0 \delta_3 \delta_3(e^h) b_0 
-2r^3b_0 \delta_3 \delta_3(e^h) b_0 b_0\\
+4r^9e^h b_0 b_0 b_0 \delta_4(e^h) e^h b_0 b_0 \delta_4(e^h) b_0 
+2r^9 e^h b_0 b_0 e^h b_0 b_0 \delta_4(e^h) b_0 \delta_4(e^h) b_0  \\
+4r^9 e^h b_0 b_0 \delta_4(e^h) e^h b_0 b_0 b_0 \delta_4(e^h) b_0  
+2r^9 e^h b_0 b_0 \delta_4(e^h) e^h b_0 b_0 \delta_4(e^h) b_0 b_0  \\
+4r^9 e^h b_0 e^h b_0 b_0 b_0 \delta_4(e^h) b_0 \delta_4(e^h) b_0  
+2r^9 e^h b_0 e^h b_0 b_0 \delta_4(e^h) b_0 b_0 \delta_4(e^h) b_0  \\
+2r^9 e^h b_0 e^h b_0 b_0 \delta_4(e^h) b_0 \delta_4(e^h) b_0 b_0  
+6r^9 e^h e^h b_0 b_0 b_0 b_0 \delta_4(e^h) b_0 \delta_4(e^h) b_0  \\
+2r^9 e^h e^h b_0 b_0 b_0 \delta_4(e^h) b_0 b_0 \delta_4(e^h) b_0  
+2r^9 e^h e^h b_0 b_0 b_0 \delta_4(e^h) b_0 \delta_4(e^h) b_0 b_0  \\
-r^7 b_0 b_0 \delta_4(e^h)e^h b_0 b_0 \delta_4(e^h) b_0 
-2r^7 b_0 \delta_4(e^h) e^h b_0 b_0 b_0 \delta_4(e^h) b_0  \\
-r^7 b_0 \delta_4(e^h) e^h b_0 b_0 \delta_4(e^h) b_0 b_0  
-16r^7 e^h b_0 b_0 b_0 \delta_4(e^h) b_0 \delta_4(e^h) b_0  \\
-r^7 e^h b_0 b_0 e^h b_0 b_0 \delta_4 \delta_4(e^h) b_0  
-8r^7 e^h b_0 b_0 \delta_4(e^h) b_0 b_0 \delta_4(e^h) b_0  \\
-8r^7 e^h b_0 b_0 \delta_4(e^h) b_0 \delta_4(e^h) b_0 b_0  
-2r^7 e^h b_0 e^h b_0 b_0 b_0 \delta_4 \delta_4(e^h) b_0  \\
-r^7 e^h b_0 e^h b_0 b_0 \delta_4 \delta_4(e^h) b_0 b_0  
-3r^7 e^h e^h b_0 b_0 b_0 b_0 \delta_4 \delta_4(e^h) b_0  \\
-r^7 e^h e^h b_0 b_0 b_0 \delta_4 \delta_4(e^h) b_0 b_0  
+5/2r^5 b_0b_0 \delta_4(e^h) b_0 \delta_4(e^h) b_0  \\
+2r^3 b_0 b_0 \delta_4(e^h) e^{-h} \delta_4(e^h) b_0  
+5/2r^5 b_0 \delta_4(e^h) b_0 b_0 \delta_4(e^h) b_0 \\ 
+5/2r^5 b_0 \delta_4(e^h) b_0 \delta_4(e^h) b_0 b_0  
+2r^3 b_0 \delta_4(e^h) e^{-h} \delta_4(e^h) b_0 b_0 \\ 
+6r^5 e^h b_0 b_0 b_0 \delta_4 \delta_4(e^h) b_0  
+3r^5 e^h b_0 b_0 \delta_4 \delta_4(e^h) b_0 b_0  \\
-2r^3 b_0 b_0 \delta_4 \delta_4(e^h) b_0 
-2r^3b_0 \delta_4 \delta_4(e^h) b_0 b_0.
$\\

We use the rearrangement lemma of \cite{conmos2} to integrate the 
above terms with respect to $r$ from $0$ to $\infty$. We recall this 
lemma here:

\begin{lemma}

For any  $m=(m_0,m_1, \dots, m_\ell) \in \mathbb{Z}^{\ell+1}_{>0}$ 
and elements $\rho_1, \dots, \rho_\ell \in C^\infty(\mathbb{T}_\theta^4)$,  
we have 
\[
\int_0^\infty \frac{u^{|m| -2}}{  (e^hu+1)^{m_0}} \prod_{1}^\ell 
\rho_j (e^hu+1)^{-m_j} \,du 
= e^{-(|m| -1)h} F_m (\Delta_{(1)}, \dots, \Delta_{(\ell)}) \Big (
\prod_1^\ell \rho_j \Big), 
\]
where
\[
 F_m(u_1, \dots, u_\ell) 
=\int_0^\infty \frac{x^{|m| -2}}{(x+1)^{m_0}}  \prod_1^l \Big (
x \prod_{1}^j u_h+1 \Big)^{-m_j}\, dx, 
\]
\[
\Delta(a)=e^{-h}ae^h, \qquad a \in C(\mathbb{T}_\theta^4).
\]
Here, $\Delta_{(j)}$ signifies the action of $\Delta$ on $\rho_j$.

\end{lemma}

Then, we use the identities (cf. \cite{conmos2, fatkha2})
\begin{eqnarray}
&&e^{-h}\delta_i(e^h) = g_1(\Delta) (\delta_i(h)),  \nonumber \\
&&e^{-h} \delta_i^2(e^h) =  g_1(\Delta) (\delta_i^2(h)) +  
2g_2(\Delta_{(1)}, \Delta_{(2)})(\delta_i(h) \delta_i(h)), \nonumber
\end{eqnarray}
where
\begin{equation} \label{gfunction}
g_1(u)=\frac{u-1}{\log u}, \qquad
g_2(u, v)= \frac{u (v-1) \log (u)-(u-1) \log (v)}{\log (u) \log (v) 
(\log (u)+\log(v))}, \nonumber
\end{equation}
and obtain the final formula for the scalar curvature in terms of 
$\nabla=\log \Delta$ and $h$, which is recorded 
in the following theorem.

\begin{theorem} \label{SCthm}

The scalar curvature $R$ of $\mathbb{T}_\theta^4$, up to a factor 
of $\pi^2$, is equal to  

\begin{eqnarray} \label{SC} 
e^{-h} K(\nabla) \Big (\sum_{i=1}^4 \delta_i^2(h) \Big ) +e^{-h} 
H(\nabla_{(1)}, \nabla_{(2)})\Big (\sum_{i=1}^4 \delta_i(h)^2 \Big ), 
\end{eqnarray}
where 
\begin{eqnarray}  
&&K(s)= \frac{1-e^{-s}}{2 s}, \nonumber \\
&&H(s,t)=-\frac{e^{-s-t} \left(\left(-e^s-3\right) s \left(e^t-1\right)+
\left(e^s-1\right) \left(3 e^t+1\right) t\right)}{4 s t (s+t)}.  \nonumber 
\end{eqnarray}

\end{theorem}

We analyse the functions $K, H,$ which describe the scalar 
curvature of $\mathbb{T}_\theta^4$, as follows.  The Taylor 
expansion of $K$ at $0$ is of the form 
\[
K(s)=\frac{1}{2}-\frac{s}{4}+\frac{s^2}{12}-
\frac{s^3}{48}+\frac{s^4}{240}-\frac{s^5}{1440}+
O\left(s^6\right).
\]
We have $\lim_{s \to -\infty} K(s)=\infty, \lim_{s \to \infty} K(s)=0$, and 
here is the graph of this function:

\vskip 0.5 cm
%\begin{figure}
\begin{center}
\includegraphics[scale=0.7]{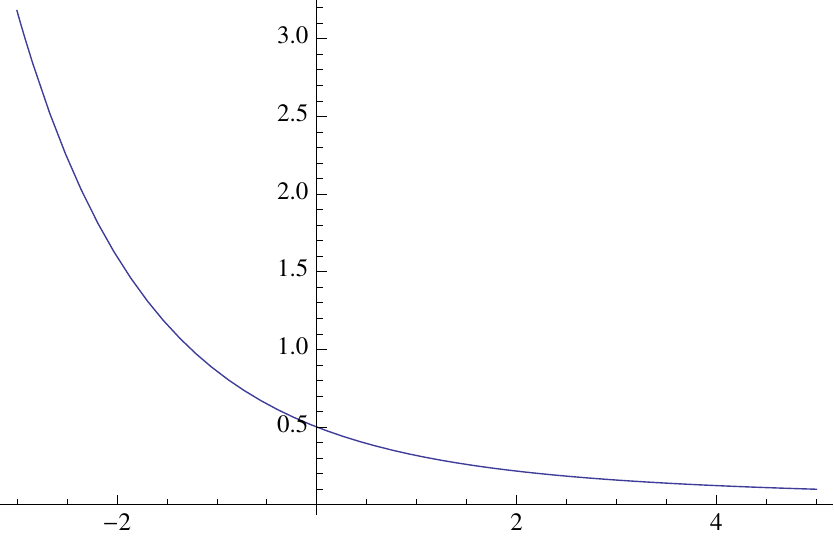}\\{Graph of the function $K$.}
%\caption{Graph of $k$.}
\label{GofK}
\end{center}
%\end{figure}
\vskip 0.5 cm

The function $H$ has the following Taylor expansion at $(0, 0)$:  
\begin{eqnarray} 
H(s,t)&=&
\left(-\frac{1}{4}+\frac{t}{24}-\frac{t^3}{480}+O\left(t^4\right)\right)+s
\left(\frac{5}{24}-\frac{t}{16}+\frac{t^2}{80}-\frac{t^3}{576}+O\left(
t^4\right)\right)\nonumber \\
&&+s^2 \left(-\frac{1}{12}+\frac{7t}{240}-\frac{t^2}{144}+\frac{5
t^3}{4032}+O\left(t^4\right)\right) \nonumber \\
&&+s^3 \left(\frac{11}{480}-\frac{5
t}{576}+\frac{t^2}{448}-\frac{t^3}{2304}+O\left(t^4\right)\right)+
O \left(s^4\right). 
\nonumber
\end{eqnarray}
Here is the graph of $H$ in a neighborhood of the origin:

\vskip 0.5 cm
%\begin{figure}
\begin{center}
\includegraphics[scale=0.5]{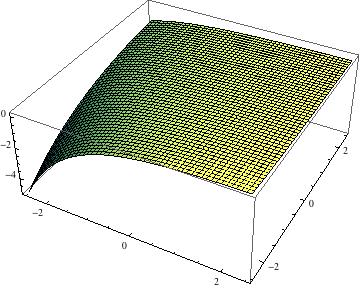}\\{Graph of the function $H$.}
%\caption{Graph of $H$.}
\label{GofH}
\end{center}
%\end{figure}
\vskip 0.5 cm

\noindent
This function is not bounded from below on the main diagonal 
as for 
\[
H(s, s)=-\frac{e^{-2 s} \left(e^s-1\right)^2}{4 s^2},
\]
we have $\lim_{s \to -\infty} H(s, s)= - \infty$, and 
$\lim_{s\to \infty} H(s, s)=0$. At $0$ we have 
the Taylor expansion 
\[
H(s, s)=-\frac{1}{4}+\frac{s}{4}-\frac{7 s^2}{48}+\frac{s^3}{16}-\frac{31
   s^4}{1440}+\frac{s^5}{160}+O\left(s^6\right), 
\]
and here is the graph of this function:

\vskip 0.5 cm
%\begin{figure}
\begin{center}
\includegraphics[scale=0.7]{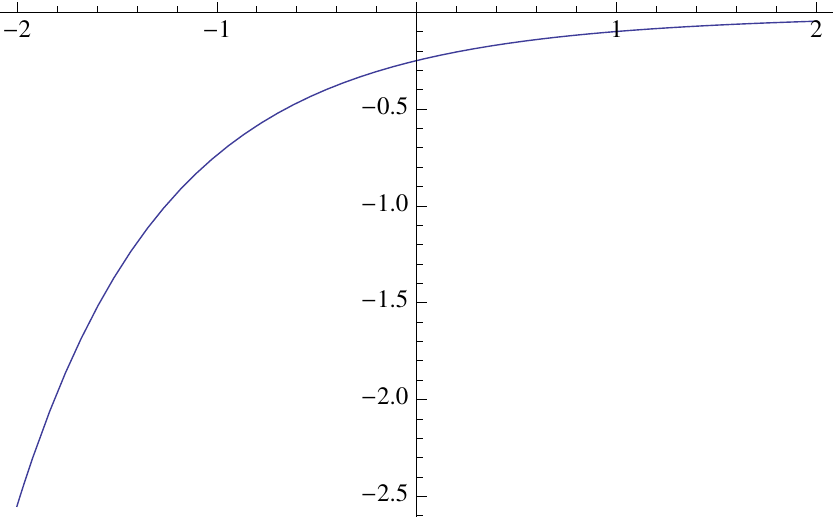}\\{Graph of the map $s \mapsto H(s, s)$.}
%\caption{Graph of $H(s, s)$.}
\end{center}
%\end{figure}
\vskip 0.5 cm

On the other diagonal, $H$ is neither bounded below nor bounded 
above as we have 
\[ 
H(s, -s)= \frac{-4 s-3 e^{-s}+e^s+2}{4 s^2}, 
\]
which implies that 
$\lim_{s \to -\infty} H(s, -s) = -\infty, \lim_{s \to \infty} H(s, -s) = \infty.$ This 
function has the following Taylor expansion at $0$
\[
H(s,-s)=-\frac{1}{4}+\frac{s}{6}-\frac{s^2}{48}+
\frac{s^3}{120}-\frac{s^4}{1440}+\frac{s^5}{5040}+O\left(s^6\right), 
\]
and here is its graph:

\vskip 0.5 cm
%\begin{figure}
\begin{center}
\includegraphics[scale=0.7]{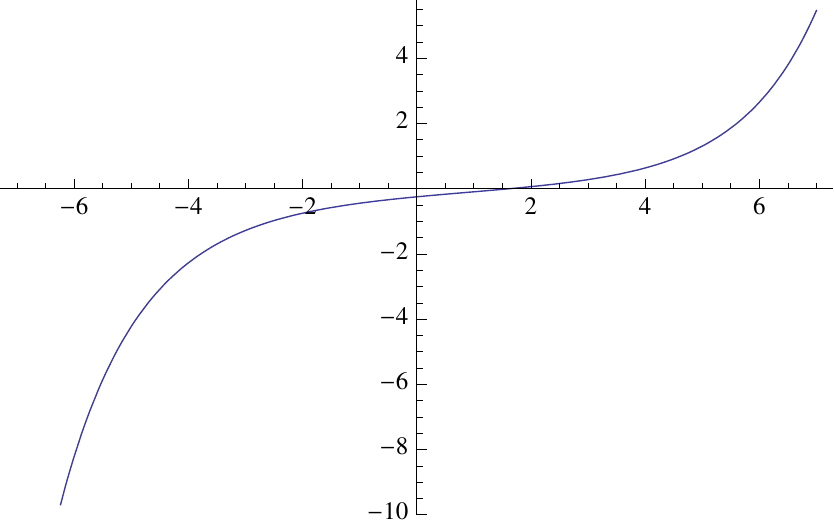}\\{Graph of the map  $s \mapsto H(s, -s)$.}
%\caption{Graph of  $H(s, -s)$.}
\label{1VarFunc}
\end{center} 
%\end{figure}
\vskip 0.5 cm

\begin{remark}

Since $K(0)=1/2$ and $H(0, 0)=-1/4$, in the commutative 
case, the scalar curvature given by \eqref{SC} reduces to 
\[
\frac{\pi^2}{2} \sum_{i=1}^4 \big (  \delta_i^2(h) - 
\frac{1}{2} \delta_i(h)^2 \big ),  
\]
which, up to a normalization factor, is the scalar curvature 
of the ordinary $4$-torus equipped with the metric 
$e^{-h}(dx_1^2+\cdots+dx_4^2)$. 

\end{remark}

\subsection{Einstein-Hilbert action for $\mathbb{T}_\theta^4$.}

A natural analogue of the Einstein-Hilbert action for $\mathbb{T}_\theta^4$ 
is $\varphi_0(R)$, where $R$ is the scalar curvature given by \eqref{SC}. In 
the following theorem we find an explicit formula for this action.

\begin{theorem}

A local expression for the Einstein-Hilbert action for 
$\mathbb{T}_\theta^4$, up to a factor of $\pi^2$, is given by 
\begin{equation} \label{EinsHilb}
\varphi_0(R) = \frac{1}{2} \varphi_0 \Big ( 
\sum_{i=1}^4 e^{-h} \delta_i^2(h) \Big )  +  \varphi_0 \Big ( 
\sum_{i=1}^4 G(\nabla)(e^{-h} \delta_i(h) ) \delta_i(h)  \Big ),
\end{equation}
where  
\[ 
G(s) =\frac{-4 s-3 e^{-s}+e^s+2}{4 s^2}.   
\]

\begin{proof}

Let us recall from Theorem \ref{SCthm} that up to a factor of $\pi^2$  
\[
R=e^{-h} K(\nabla) \Big (\sum_{i=1}^4 \delta_i^2(h) \Big ) +e^{-h} 
H(\nabla_{(1)}, \nabla_{(2)})\Big (\sum_{i=1}^4 \delta_i(h)^2 \Big ).
\]
Writing $K$ as a Fourier transform 
\[
K(s) = \int e^{- i u s } f(u) \, du,
\]
we have  
\[
\varphi_0 \big (e^{-h} K(\nabla) (\delta_i^2(h)) \big ) = 
\int \varphi_0 \big (   e^{-h} \Delta^{-iu} 
(\delta_i^2(h))  \big ) f(u)\, du 
= K(0) \, \varphi_0(e^{-h} \delta_i^2(h)).
\]
Also, by writing $H$ as a Fourier transform
\[
H(s,t) = \int e^{-i(su+tv)} g(u, v) \,du\,dv,
\]
we have
\begin{eqnarray}
&& \varphi_0\big ( e^{-h} H(\nabla_{(1)}, \nabla_{(2)})( \delta_i(h)^2) \big ) \nonumber \\
&&=  \int \varphi_0\big(  e^{-h} \Delta^{-iu}( \delta_i(h)) \Delta^{-iv} (\delta_i(h)) \big )
g(u, v) \, du \, dv\nonumber \\
&&= \int \varphi_0 \big ( \Delta^{-i(u-v)}(e^{-h} \delta_i(h)) \delta_i(h) \big ) g(u, v)
\, du \, dv \nonumber \\
&&= \varphi_0 \big (  H(\nabla, -\nabla) (e^{-h} \delta_i(h)) \delta_i(h)   
\big ). \nonumber
\end{eqnarray}
Therefore
\begin{eqnarray}
\varphi_0(R) = \frac{1}{2} \varphi_0 \Big ( 
\sum_{i=1}^4 e^{-h} \delta_i^2(h) \Big )  +  \varphi_0 \Big ( 
\sum_{i=1}^4 G(\nabla)(e^{-h} \delta_i(h) ) \delta_i(h)  \Big ), \nonumber
\end{eqnarray}
where
\[ 
G(s) =H(s, -s) = \frac{-4 s-3 e^{-s}+e^s+2}{4 s^2}.   
\]

\end{proof}

\end{theorem}

\subsection{Extremum of the Einstein-Hilbert action for $\mathbb{T}_\theta^4$.}

We show that the Einstein-Hilbert action $\varphi_0(R)$ attains its 
maximum if and only if the Weyl factor $e^{-h}$ is a constant. This is 
done by combining the two terms in the explicit formula \eqref{EinsHilb} for 
$\varphi_0(R)$, and observing that it can be expressed by a 
non-negative function. We note that the function $G$ in \eqref{EinsHilb}, 
which was analysed in Subsection \ref{SCsub},  is neither bounded below 
nor bounded above.

\begin{theorem}

The maximum of the Einstein-Hilbert action is equal to $0$, and it is 
attained if and only if the Weyl  factor is a constant. That is, 
for any Weyl  factor $e^{-h}, h=h^*\in C^\infty(\mathbb{T}_\theta^4),$ 
we have 
\[
\varphi_0(R) \leq 0,
\]
and the equality happens if and only if $h$ is a constant.

\begin{proof}

We can combine the two terms in \eqref{EinsHilb}  as follows. We have

\begin{eqnarray}
\varphi_0(e^{-h} \delta_i^2(h)) &=& 
- \varphi_0(\delta_i(e^{-h}) \delta_i(h)) \nonumber \\
&=&\varphi_0(e^{-h}\delta_i(e^{h}) e^{-h} \delta_i(h)) \nonumber \\ 
&=& \varphi_0( e^{-h} \delta_i(h) e^{-h}\delta_i(e^{h}) )  \nonumber \\ 
&=& \varphi_0 \Big (e^{-h}\delta_i(h)  \frac{\Delta-1}{\log \Delta}(\delta_i(h)) \Big ) \nonumber \\ 
&=& \varphi_0 \Big (e^{-h}\delta_i(h)  \frac{e^\nabla -1}{\nabla}(\delta_i(h)) \Big) \nonumber \\ 
&=& \varphi_0 \Big (e^{-h}\frac{e^{-\nabla} -1}{-\nabla}(\delta_i(h))  \delta_i(h) \Big) \nonumber. 
\end{eqnarray}
The last equality follows from the fact that for any entire function 
F, one has 
\[
\varphi_0(e^{-h} a F(\nabla)(b)) =\varphi_0(e^{-h} F(-\nabla)(a) b), 
\qquad  a, b \in C(\mathbb{T}_\theta^4).  
\]
Therefore we can write \eqref{EinsHilb} as 
\begin{eqnarray} \label{combined}
\varphi_0(R)&=&   
\sum_{i=1}^4 \varphi_0 \Big(\frac{1}{2}  e^{-h}\frac{e^{-\nabla} -1}{-\nabla}(\delta_i(h))  
\delta_i(h)+G(\nabla)(e^{-h} \delta_i(h) ) \delta_i(h) \Big ) \nonumber \\
&=& \sum_{i=1}^4  \varphi_0 \big ( e^{-h} T(\nabla)(\delta_i(h) ) 
\delta_i(h) \big ), 
\end{eqnarray}
where
\[
T(s)=\frac{1}{2} \frac{e^{-s}-1}{-s}+G(s)=\frac{-2 s+e^s-e^{-s} (2 s+3)+2}{4 s^2}.
\]

The Taylor expansion of $T$ at 0 is of the form 
\[
T(s)=\frac{1}{4}-\frac{s}{12}+\frac{s^2}{16}-\frac{s^3}{80}+
\frac{s^4}{288}-\frac{s^5}{20
   16}+O\left(s^6\right).
\]
We have $\lim_{s \to \infty}T(s)=\infty, \lim_{s \to - \infty} = \infty$. 
Moreover, this function is non-negative as its absolute minimum is approximately 
$0.218207$ which is attained around $s = 0.812394$. Here is the graph of this 
function:

\vskip 0.5 cm
%\begin{figure}
\begin{center}
\includegraphics[scale=0.7]{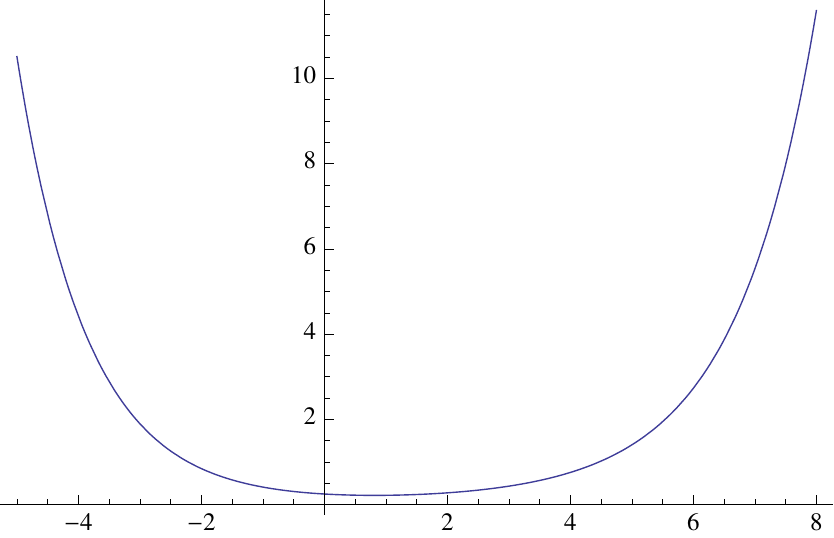} \\{Graph of the function $T$.}
%\caption{Graph of $T1$.}
%\label{1VarFunc}
\end{center}
%\end{figure}
\vskip 0.5 cm

Since $T$ is a non-negative function and $\delta_i(h)^*=-\delta_i(h)$ for 
$i=1, \dots,4$, we have 
\[
\varphi_0 ( e^{-h} T(\nabla)(\delta_i(h)) \delta_i(h))= 
- \varphi_0( e^{-h} T(\nabla)(\delta_i(h)) \delta_i(h)^*) \leq 0. 
\]
In the last inequality we have used the fact that $\nabla$ is a 
selfadjoint operator with respect to the inner product 
\[
(a, b) =  \varphi_0(e^{-h} b^*a ), \qquad a,b \in C(\mathbb{T}_\theta^4).
\]
Also, faithfulness of $\varphi_0$ implies that for each 
$i=1, \dots, 4,$ 
\[ 
\varphi_0( e^{-h} T(\nabla)(\delta_i(h)) \delta_i(h))=0 
\]
if and only if $\delta_i(h)=0$. Thus, it follows from \eqref{combined} that 
\[ 
\varphi_0(R)= \sum_{i=1}^4  \varphi_0 \big ( e^{-h} 
T(\nabla)(\delta_i(h)) \delta_i(h) \big ) \leq 0,
\] 
and the equality happens if and only $\delta_i(h)=0$ for all $i=1, 
\dots, 4,$ which holds if and only if $h$ is a constant. 

\end{proof}

\end{theorem}

$ $\\
\noindent
Department of Mathematics, Western University \\  
London, ON, Canada, N6A 5B7 \\
{\it E-mail address}:  ffathiz@uwo.ca \\

\noindent
Department of Mathematics, Western University \\  
London, ON, Canada, N6A 5B7 \\
{\it E-mail address}:  masoud@uwo.ca

\end{document}